\documentclass[11pt]{amsart}
\usepackage{amsmath,amssymb,latexsym,dsfont}
\usepackage[left=2.6cm,right=2.6cm,top=2.7cm,bottom=2.7cm]{geometry}

\usepackage{amsmath,amsfonts,amssymb,epsfig, textcomp}
\usepackage{graphicx,tikz}
\usepackage{enumitem}
\usepackage{hyperref, cleveref} 
\usepackage{cases,listings}
\usepackage{color}
\usepackage{mathabx,esint}

\newtheorem{theorem}{Theorem}[section]
\newtheorem{lemma}{Lemma}[section]

\newtheorem{proposition}{Proposition}[section]
\newtheorem{corollary}{Corollary}[section]
\newtheorem{remark}{Remark}[section]

\numberwithin{equation}{section}

      \newcommand{\hu}{\hat u}

   \newcommand{\ty}{\widetilde y}

      \newcommand{\N}{\mathbb{N}}

      \newcommand{\mC}{\mathbb{C}}

      \makeatletter
      \def\@setcopyright{}
      \def\serieslogo@{}
      \makeatother

\newcommand{\cL}{\mathcal L}

\newcommand{\tu}{\widetilde u}

\newcommand{\be}{\begin{equation}}
\newcommand{\ee}{\end{equation}}

 \newcommand{\hx}{\hat x}

 \newcommand{\tx}{\widetilde{x}}

\newcommand{\cD}{{\mathcal D}}

\newcommand{\mH}{\mathbb{H}}

\newcommand{\mU}{\mathbb{U}}

\newcommand{\mX}{\mathbb{X}}

\title[Turnpike property of linear control systems]{Turnpike property of linear quadratic control problems with unbounded control operators}

\author[H.-M. Nguyen]{Hoai-Minh Nguyen}
\address[H.-M. Nguyen]{Sorbonne Université, Université Paris Cité, \newline\indent
CNRS, INRIA, Laboratoire Jacques-Louis Lions, LJLL, 
\newline\indent
F-75005 Paris, France}
\email{hoai-minh.nguyen@sorbonne-universite.fr}

\author[E. Tr\'elat]{Emmanuel Tr\'elat}
\address[E. Tr\'elat]{Sorbonne Université, Université Paris Cité, \newline\indent
CNRS, INRIA, Laboratoire Jacques-Louis Lions, LJLL, 
\newline\indent
F-75005 Paris, France}
\email{emmanuel.trelat@sorbonne-universite.fr}

\begin{document} 

\begin{abstract} 
We establish the turnpike property for linear quadratic control problems for which the control operator is admissible and may be unbounded, under quite general and natural assumptions.  
The turnpike property has been well studied for bounded control operators, based on the theory of differential and algebraic Riccati equations. For unbounded control operators, there are only few results, limited to some special cases of hyperbolic systems in dimension one or to analytic semigroups. 
Our analysis is inspired by the pioneering work of Porretta and Zuazua \cite{PZ13}. 
We start by approximating the admissible control operator with a sequence of bounded ones. We then prove the convergence of the approximate problems to the initial one in a suitable sense. Establishing this convergence is the core of the paper. It requires to revisit in some sense the linear quadratic optimal control theory with admissible control operators, in which the roles of energy and adjoint states, and the connection between infinite-horizon and finite-horizon optimal control problems with an appropriate final cost are investigated.   
\end{abstract}

\maketitle

\tableofcontents

\section{Introduction}

\subsection{Setting and main results}

Let $\mH$ and $\mU$ be two Hilbert spaces standing for the state space and the control space, respectively. 
The corresponding scalar products are  $\langle \cdot, \cdot \rangle_{\mH}$ and $\langle \cdot, \cdot \rangle_{\mU}$, and the corresponding norms are $\| \cdot \|_{\mH}$ and $\| \cdot \|_{\mU}$. 
Throughout the paper, the dual $\mH'$ of $\mH$ (resp., the dual $\mU'$ of $\mU$) is identified with $\mH$) (resp. with $\mU$).
Given any two Hilbert spaces $\mX_1$ and $\mX_2$, we denote by $\cL(\mX_1, \mX_2)$ the Banach space of all bounded linear applications from $\mX_1$ to $\mX_2$ endowed with the usual operator norm, and we denote $\cL(\mX_1, \mX_1)$ by $\cL(\mX_1)$.

Let $A:\cD(A)\rightarrow\mH$ be a linear operator generating on $\mH$ a strongly continuous semigroup $( e^{tA} )_{t \ge 0} \subset \cL(\mH)$ (see \cite{EN00,Pazy}).
The adjoint operator $A^*:\cD(A^*)\rightarrow\mH$ generates the adjoint semigroup $( e^{tA^*} )_{t \ge 0} \subset \cL(\mH)$, with $e^{tA^*}=(e^{tA})^*$ for every $t\ge 0$. The domain $\cD(A^*)$, equipped with the scalar product 
$$
\langle z_1, z_2 \rangle_{\cD(A^*)} = \langle z_1, z_2 \rangle_{\mH} +  \langle A^* z_1, A^* z_2 \rangle_{\mH} \qquad \forall z_1, z_2 \in \cD(A^*),
$$
is a Hilbert space. The dual $\cD(A^*)'$ of $\cD(A^*)$ with respect to the pivot space $\mH$ satisfies
$\cD(A^*) \subset \mH \subset \cD(A^*)'$ with continuous and dense inclusions.
 
Given any $x_0\in\mH$ and given a control operator
$$
B \in \cL(\mU, \cD(A^*)'),
$$
we consider the control system
\be \label{CS}
\begin{split}
x'(t) &= Ax(t) + Bu(t), \qquad t>0, \\
x(0) &= x_0, 
\end{split} 
\ee
where, at time $t$, the control is $u(t) \in \mU$ and the state is $x(t) \in \mH$ (see \cite{BDDM07, Coron07, CZ95, Lions71, Trelat_2024, TW09, Zabczyk08} for references on control theory in infinite dimension). Interesting aspects of controllability and stabilization properties of \eqref{CS} can be found in \cite{Coron07,CZ95, EN00, Lions71,Slemrod74,TWX20,TW09,WZ98,Zabczyk08} (see also references therein). 

We assume that $B$ is an {\it admissible} control operator with respect to the semigroup $( e^{tA})_{t \ge 0}$ (see \cite{TW09}) in the sense that, for all $u \in L^2((0, T); \mU)$, it holds that 
$$
\varphi \in C([0, T]; \mH) \mbox{ where } \varphi(t) = \int_0^t e^{(t - s)A} B u(s) \, ds. 
$$
As a consequence of the closed graph theorem (see, e.g., \cite{Brezis-FA}), one has 
$$
\| \varphi\|_{C([0, T]; \mH)} \le c_T \| u \|_{L^2((0, T); \mU)},  
$$
for some positive constant $c_T$ not depending on $\varphi$.

Let now
$$
C \in \cL(\mH),  
$$
be an observation operator. Let $z \in \mH$ be fixed.

On the one part, given any $T>0$, we consider the \emph{linear quadratic optimal control problem} (dynamical optimization problem)
\be \label{optimal-control-T}
\inf_{u \in L^2((0, T); \mU)} \int_0^T \big( \|Cx(t) - z\|_{\mH}^2 + \|u(t)\|_{\mU}^2  \big) \, dt,  
\ee
where $x \in C([0, T]; \mH)$ is the (unique) solution of \eqref{CS}.
By strict convexity, there exists a unique optimal control $u_{T,\mathrm{opt}} \in L^2((0, T); \mU)$ solution of \eqref{optimal-control-T}. Let $x_{T,\mathrm{opt}} \in C([0, T]; \mH)$ be the corresponding solution and denote $y_{T,\mathrm{opt}} $ the corresponding adjoint state, i.e., $y_{T,\mathrm{opt}} \in  C([0, T]; \mH)$ is the solution of the system 
\begin{equation*}
\begin{split}
y_{T,\mathrm{opt}}'(t) &= - A^* y_{T,\mathrm{opt}}(t) - C^* (C x_{T,\mathrm{opt}}(t) - z) \quad\mbox{ on } (0, T), \\
y_{T,\mathrm{opt}}(T) &= 0 
\end{split} 
\end{equation*}
(see, e.g., \Cref{pro-Opt-T} for the role of the adjoint state $y_{T,\mathrm{opt}}$). 

On the other part, we consider the \emph{stationary optimization problem}
\be \label{optimal-control-stationary}
\inf_{(x, u) \in \mH \times \mU \atop A x + B u =0} \left( \|Cx - z\|_{\mH}^2 + \|u\|_{\mU}^2  \right) .
\ee
Under appropriate assumptions (see \Cref{pro-stationary} below), this problem has a unique optimal solution $(\bar x, \bar u) \in \mH \times \mU$ and a unique Lagrange multiplier $\bar y \in \cD(A^*)$.

The objective of this paper is to establish the \emph{exponential turnpike property}, namely that, when $T$ is large, the optimal solution of \eqref{optimal-control-T} remains exponentially close to the optimal solution of the stationary problem \eqref{optimal-control-stationary}, except at the extremities of the time interval $[0,T]$.
The turnpike property for \eqref{optimal-control-T} has already been established in the literature (see further) but mainly for bounded control operators $B$, i.e., when $B \in L(\mU, \mH)$, or under strong assumptions on the semigroup. Our main goal is to establish the exponential turnpike property under quite general and natural assumptions on the triple $(A,B,C)$ for which the boundedness of  $B$ is not required (see Theorem \ref{thm-main} below).

\medskip
To state the turnpike property, we first discuss
the existence and uniqueness of the optimal solution and of its Lagrange multiplier for the stationary problem \eqref{optimal-control-stationary}.  

\begin{proposition} \label{pro-stationary} Let $z \in \mH$. 
We assume that 
\be \label{pro-stationary-AC}
\ker A \cap \ker C = \{ 0 \}
\ee
and that there exists $T_0>0$ such that $(A, C)$ is finite-time observable in time $T_0$, i.e., that there exists $c >0$ such that
$$
\int_{0}^{T_0} \|C e^{t A} \xi \|_{\mH}^2  \, dt \ge c \|e^{T_0 A} \xi\|_{\mH}^2 \qquad\forall \xi \in \mH.
$$
Then, there exists a unique pair $(\bar x, \bar u) \in \mH \times \mU$ such that  
$A \bar x + B \bar u = 0$\footnote{This means that $\langle \bar x, A^*  \varphi \rangle_{\mH} + \langle \bar u, B^* \varphi \rangle_{\mU} = 0$ for every $\varphi \in \cD(A^*)$.} 
and
\be \label{pro-stationary-cl1}
\|C \bar x - z\|_{\mH}^2 + \|\bar u\|_{\mU}^2 = \inf_{(x, u) \in \mH \times \mU \atop A x + B u =0} \left( \|Cx - z\|_{\mH}^2 + \|u\|_{\mU}^2 \right).
\ee
If we assume in addition that
$$
\ker A^* \cap \ker B^* = \{0\}
$$
and that the pair $(A^*, B^*)$ is finite-time observable in time $T_0$, i.e., that there exists $c >0$ such that
\be \label{pro-stationary-A*B*}
\int_{0}^{T_0} \|B^* e^{t A^*} \xi \|_{\mU}^2  \, dt \ge c \|e^{T_0 A} \xi\|_{\mH}^2 \qquad\forall \xi \in \mH, 
\ee
then there exists a unique $\bar y \in \cD(A^*)$ such that 
\be \label{pro-stationary-cl2}
- A^* \bar y - C^* (C \bar x - z) = 0 \quad \mbox{ and } \quad \bar u = - B^* \bar y. 
\ee
\end{proposition}

The  proof of \Cref{pro-stationary} is given in \Cref{sect-stationary}.

\begin{remark} \rm  As seen in the proof of \Cref{pro-stationary}, the observability of $(A, C)$ implies the existence of $(\bar x, \bar u) \in \mH \times \mU$ satisfying \eqref{pro-stationary-cl1}. The assumption $\ker A \cap \ker C = \{ 0 \}$ implies the uniqueness of such a pair $(\bar x, \bar u)$. The additional assumption of observability of $(A^*, B^*)$ implies the existence of $\bar y \in \mH$ satisfying \eqref{pro-stationary-cl2}, and  the assumption $\ker A^* \cap \ker B^* = \{ 0 \}$ implies the uniqueness of $\bar y $. 
\end{remark}

\begin{remark} \rm The unique $\bar y \in \mH$ given in \Cref{pro-stationary} is the Lagrange multiplier of the constrained optimization problem \eqref{optimal-control-stationary}.
A general theory on this topic can be found in, e.g., \cite{Luenberger69}, under the assumption that the operator $\mH \times \mU\ni(x,u)\mapsto A x + Bu$ be surjective.
\end{remark}

Our main result, proved in \Cref{sect-main}, is the following. 

\begin{theorem} \label{thm-main} 
We make the following assumptions:
\begin{enumerate}[label=$\bf (H_\arabic*)$]
\item\label{H2} The pairs $(A, C)$ and $(A^*, B^*)$ are finite-time observable in some time $T_0>0$,  and 
\be \label{thm-main-ker}
\ker A \cap \ker C = \{ 0 \} \quad \mbox{ and } \quad \ker A^* \cap \ker B^* = \{ 0 \}. 
\ee
\item\label{H3} There exists $\delta>0$ such that
$$
C^* C \ge \delta\, \mathrm{id}, \quad\mbox{i.e., } \| C \eta \|_{\mH}^2 \ge \delta \| \eta \|_{\mH}^2 \quad\forall\eta \in \mH. 
$$
\end{enumerate}
Then, we have the exponential turnpike property: there exist $c>0$ and $\lambda>0$ (depending only on $A$, $B$, $C$) such that, for any $T > T_0$, and any $x_0,z \in \mH$,
\begin{multline} \label{thm-main-cl}
\|x_{T,\mathrm{opt}}(t) - \bar x\|_{\mH} + \|y_{T,\mathrm{opt}}(t) - \bar y\|_{\mH} + \|u_{T,\mathrm{opt}} - \bar u\|_{L^2(I_t ; \mU)}   \\
\le c  \big(\Vert x_0-\bar x\Vert_\mH + \Vert\bar y\Vert_\mH \big)  \big( e^{-\lambda t} + e^{- \lambda(T - t)} \big) 
\qquad\forall t\in[0, T] ,
\end{multline}
where $u_{T,\mathrm{opt}}$ is the solution of the optimal control problem of \eqref{optimal-control-T}, $x_{T,\mathrm{opt}}$ is the corresponding solution of \eqref{CS}, $y_{T,\mathrm{opt}}$ is the corresponding adjoint state, $(\bar x, \bar u, \bar y)$ is the (unique) triple solution of the stationary problem \eqref{optimal-control-stationary}, given by \Cref{pro-stationary}, and $I_t = (t,T-t)$ if $t\leq T/2$ and $I_t = (T-t,t)$ if $t\geq T/2$.
\end{theorem}

\begin{remark} \rm As seen in the proof of the theorem, the decay rate $\lambda$ can be chosen as the the decay rate of the semigroup $S_{\infty,\mathrm{opt}}(t)$ defined in \Cref{pro-def-Sopt} further. This decay rate is optimal in the exponential turnpike property \eqref{thm-main-cl} (see \Cref{pro-technical}). Assumption \ref{H3} is used to ensure that $S_{\infty,\mathrm{opt}}(t)$ exponentially decays, see \Cref{pro-def-Sopt}. Any assumption ensuring this can replace Assumption \ref{H3} in \Cref{thm-main}. 
\end{remark}

\begin{remark} \rm It is worth noting (see \Cref{pro-stationary} and \Cref{pro-Opt-T}) that
$$
u_{T,\mathrm{opt}} - \bar u = - B^* (y_{T,\mathrm{opt}} - \bar y) \quad\mbox{ in } L^2((0, T); \mU). 
$$
When $B$ is bounded, it follows from \eqref{thm-main-cl} that 
\be \label{rem-thm-main-cl-B}
 \|u_{T,\mathrm{opt}}(t) - \bar u\|_{\mU} \le c \big(\Vert x_0-\bar x\Vert_\mH + \Vert\bar y\Vert_\mH \big) \big( e^{-\lambda t} + e^{- (T - t)} \big)  \qquad\forall t\in[0, T]. 
\ee
In contrast with the pointwise estimate \eqref{rem-thm-main-cl-B}, in the estimate \eqref{thm-main-cl} established in \Cref{thm-main} we just obtain a $L^2$ estimate on the control. This is consistent with the fact that $B$ is admissible and might be unbounded. 
\end{remark}

\begin{remark} \rm The existence and uniqueness of the optimal triple $(\bar x, \bar u, \bar y)$ for the stationary problem have been overlooked in the existing literature, even for bounded $B$. It is worth noting that the turnpike property is not meaningful unless the existence and uniqueness of such a triple are guaranteed.
\end{remark}

\subsection{Literature review}
The exponential turnpike property \eqref{thm-main-cl} has been first established in \cite{PZ13} for $B$ bounded under some additional technical assumptions. In \cite{TZZ18}, the exponential turnpike property was established for general abstract linear control systems \eqref{CS} for bounded $B$ under the assumption that $(A,B)$ is exponentially stabilizable and $(A,C)$ is exponentially detectable, or for $B$ unbounded and admissible under the assumption that the semigroup $(e^{tA})_{t \ge 0}$ is analytic. Nonlinear variants under some smallness condition are also considered there. 
The optimal decay rate was obtained in these works but the explicit dependence on $z$ and $x_0$ on the right-hand side of \eqref{thm-main-cl} was not studied there. In the periodic setting, i.e., $A$, $B$, and $C$ are periodic functions with the same period, the turnpike property was also investigated in \cite{TZZ18} for $B$ bounded. The turnpike property has been generalized to unbounded control operators in several situations. In \cite{GugatHante_2019}, the authors obtained the exponential turnpike property for $2 \times 2$ hyperbolic system in dimension one. In \cite{GugatTrelatZuazua_SCL2016}, the authors derived the turnpike property for the wave equations in dimension one with Neumann boundary controls.  
We also quote \cite{ATZ24} for results on the behavior of the cost function, related to the turnpike property. The key tool used in the works mentioned here is the theory of differential and algebraic Riccati equations for \eqref{CS} and \eqref{optimal-control-T} with $z=0$. The Riccati theory is well established and known for $B$ bounded, but in the unbounded case it is much more involved and delicate. 

In \cite{GruneSchallerSchiela_JDE2020}, the authors establish the exponential turnpike property, not by using Riccati theory but by developing a multiplier technique and a kind of integration by parts using the exponential stabilizability of $(A, B)$ and the exponential detectability of $(A, C)$. 
Most of the paper is dedicated to bounded control operators, but in \cite[Section 5]{GruneSchallerSchiela_JDE2020} the authors explain how to extend their results to the case of unbounded admissible control operators. Nevertheless, there is a problem in the proof of \cite[Lemma 6]{GruneSchallerSchiela_JDE2020}, which is one of the key ingredients of their analysis. Precisely, it is wrong that ``\cite[(19)]{GruneSchallerSchiela_JDE2020}  follows analogous by testing the adjoint equation with $\psi$ solving \cite[(16)]{GruneSchallerSchiela_JDE2020}", because it is not true that the domains $\cD(A)$ and $\cD(A+BK_B)$, where $A + B K_B$ is the operator coming from the exponential stabilizability of $(A, B)$, are the same or even have a dense intersection. The multiplier technique and the integration by parts given there thus cannot work under the assumptions made. We thank very much Manuel Schaller for a discussion regarding this issue. 

Other information on the turnpike property in infinite dimension and related methods and contributions can be found in \cite{CirantPorretta_COCV2021, Gugat_MCSS2021, FaulwasserGrune_2022, GZ22, GuglielmiLi_MCSS2024, HangZuazua_SICON2022, LanceTrelatZuazua_SCL2020,  SunYong_SICON2024, TZ15, Zaslavski06, Zaslavski15} and references therein.

\subsection{Main ideas of the proof}
Let us briefly discuss the main ideas of the proof of \Cref{thm-main}. Our analysis is strongly inspired by the work of Porretta and Zuazua \cite{PZ13}. The key idea is to establish a connection between the optimal control problem \eqref{optimal-control-T} and the optimal control problem in infinite horizon 
\be \label{optimal-control-infinite}
\inf_{u \in L^2((0, + \infty); \mU)} \int_0^\infty \left( \|Cx(t)\|_{\mH}^2 + \|u(t)\|_{\mU}^2  \right) dt ,
\ee
where $x(\cdot)$ is the solution of $x'(t) = A x(t) + Bu(t)$ on $(0, + \infty)$ such that $x(0) = x_0 \in \mH$.

Regarding this optimal control problem in infinite horizon, one has the following result (see \cite{FLT88}). 
Recall that the triple $(A, B, C)$ enjoys the finite cost condition, i.e., for every $\xi \in \mH$, there exists $u \in L^2((0, + \infty); \mU)$ such that 
$\int_0^\infty \left(\|C x(t)\|_{\mH}^2 + \|u(t) \|_{\mU}^2 \right) dt < + \infty$,
where $x(\cdot)$ is the unique solution of $x'(t) = Ax(t) + Bu(t)$ on $(0,+\infty)$ such that $x(0) = \xi$.

\begin{proposition} \label{pro-def-Sopt} 
Assume that the triple $(A, B, C)$ enjoys the finite cost condition.
For $\xi \in \mH$, let $\hu_{\infty,\mathrm{opt}, \xi}$ be the optimal control solution of the problem 
$$
\inf_{u \in L^2((0, + \infty); \mU)}  \int_0^\infty \left(\|C  x(t)\|_{\mH}^2 + \|u(t) \|_{\mU}^2 \right)  dt
\quad\mbox{ where } 
\left\{\begin{array}{rcl}
x' (t)&=& Ax(t) + Bu(t) \quad\mbox{on } (0, +\infty) , \\
x(0) &=& \xi. 
\end{array}\right.
$$
Denote $\hx_{\infty,\mathrm{opt}, \xi}$ the corresponding solution. For $t \ge 0$, define 
$$
\begin{array}{rcl}
S_{\infty,\mathrm{opt}}(t): \mH  &\to & \mH \\
 \xi & \mapsto & \hx_{\infty,\mathrm{opt}, \xi}(t). 
\end{array} 
$$
Then $(S_{\infty,\mathrm{opt}}(t) \big)_{t \ge 0} \subset \cL(\mH)$ is a strongly continuous semigroup on $\mH$. Assume in addition that either there exists $\delta>0$ such that
$$
C^* C \ge \delta\,\mathrm{id}, \quad\mbox{i.e., } \| C \eta \|_{\mH}^2 \ge \delta \| \eta \|_{\mH}^2 \quad\forall\eta\in\mH. 
$$
Then $S_{\infty,\mathrm{opt}}(t)$ decays exponentially as $t\rightarrow+\infty$.  
\end{proposition}

We next recall a property on the finite cost condition of the triple $(A, B, C)$, see \cite[Section 2]{FLT88}.

\begin{lemma} \label{lem-P} 
Assume that the triple $(A, B, C)$ enjoys the finite cost condition. Then there exists a unique symmetric $P \in \cL(\mH)$, defined by 
\be \label{lem-P-P}
\langle P \xi, \xi \rangle_{\mH} = \inf_{u \in L^2((0, + \infty); \mU)}  \int_0^\infty \left(\|C  x(t)\|_{\mH}^2 + \|u(t) \|_{\mU}^2 \right) dt, 
\ee
where $x(\cdot)$ is the solution of $x'(t) = A x(t) + Bu(t)$ on $(0, + \infty)$ such that $x(0) = \xi \in \mH$.
\end{lemma}

\begin{remark} \label{rem-FiniteCost} \rm Assume that  $(A^*,B^*)$ is observable for some positive time $T_0$. Then the control system \eqref{CS} is null-controllable for time $T_0$ (see, e.g., \cite{Coron07,TW09}) and therefore the triple $(A,B,C)$ enjoys the finite cost condition. 
\end{remark}

The following result establishes an important connection between the optimal control problem \eqref{optimal-control-T} and the optimal control problem \eqref{optimal-control-infinite}, which is instrumental for establishing \Cref{thm-main}.

\begin{proposition}\label{pro-technical} 
Under Assumption \ref{H2} of \Cref{thm-main}, let $z \in \mH$ and let $T > T_0$. Given $x_0 \in \mH$,  
set 
\be \label{pro-technical-def-hg}
h (t) = y_{T,\mathrm{opt}}(t) - \bar y  -  P (x_{T,\mathrm{opt}}(t) - \bar x)  \quad \mbox{ and } \quad g(t) = h(T - t) 
\quad\forall t\in[0, T],
\ee
where $P$ is defined by \eqref{lem-P-P}, $x_{T,\mathrm{opt}}$ is the optimal solution of \eqref{optimal-control-T},  $y_{T,\mathrm{opt}}$ is the corresponding adjoint state, and $(\bar x, \bar y)$ is given by \Cref{pro-stationary}. Then 
\be \label{pro-technical-pro-g}
g(t) = S_{\infty,\mathrm{opt}}^* (t) g(0) \qquad\forall t \ge 0,   
\ee
where $S_{\infty,\mathrm{opt}}(t)$ is defined in \Cref{pro-def-Sopt} and $S_{\infty,\mathrm{opt}}^*(t)$ is its adjoint. 
\end{proposition}

As a consequence of \Cref{pro-technical}, we derive from \Cref{pro-def-Sopt} the decay rate of $t\mapsto y_{T,\mathrm{opt}}(T- t) - \bar y  -  P (x_{T,\mathrm{opt}}(T-t) - \bar x)$. In \Cref{pro-technical}, the control operator $B$ is not required to be bounded. When $B$ is bounded, \eqref{pro-technical-pro-g} can be written as
$$
g'(t) = (A^*- P B B^*) g(t) \quad\mbox{on } (0, T)
$$
 (see \Cref{lem-technical1}), i.e., equivalently, 
\be \label{rem-tech1}
h'(t) = (-A^*+ P B B^*) h(t) \quad\mbox{on } (0, T). 
\ee
The proof of \Cref{pro-technical} is the core of our analysis. We first prove \eqref{pro-technical-pro-g} for $B$  bounded by establishing \eqref{rem-tech1}. This is inspired by the approach of Porretta and Zuazua \cite{PZ13}. We then approximate the pair $(A, B)$ in \eqref{CS} by the pair $(A, B_k)$ (i.e., we replace $B$ by $B_k$ in \eqref{CS}) and then we study the limit as $k\rightarrow+\infty$. Here $B_k$ is defined, for $k\in\N$ large enough, by
\be \label{def-Bk}
B_k = J_k B, 
\ee
where $J_k$ is the Yosida approximation of the identity with respect to $A$, defined, e.g., by $J_k x = k  \int_0^\infty e^{- k s} e^{s A} x \, ds$ for every $x \in \mH$. 

To this end, we first establish the existence and uniqueness of the optimal solutions and its Lagrange multipliers for the approximate stationary problems, and show their convergence to the optimal solution and Lagrange multiplier of the initial stationary problem (see \Cref{sect-stationary}). 

We then establish the convergence of the approximate optimal control problems to the initial control problem, which involves $S_{\infty,\mathrm{opt}}$ (see \Cref{sect-ingredient}). In this part of the analysis, to avoid dealing with the finite cost condition for the approximate control system, which might not hold, we make some modifications in comparison with \Cref{pro-technical} for the approximate systems, see \Cref{lem-technical2} in which a finite horizon optimal control problem with a suitable choice of the final cost is considered instead. The incorporation of the final cost and a connection between the optimal control problem in infinite horizon with the optimal control problem in finite horizon with an appropriate final cost (see \Cref{pro-P}) are useful in our analysis to avoid dealing with $k \to + \infty$ and $T = + \infty$ at the same time as well (see \Cref{lem-ap}).  

The introduction of $h$ and $g$ has its roots in the work \cite{PZ13}. Indeed, consider the following optimal control in finite horizon
\be \label{optimal-control-finite}
\inf_{u \in L^2((0, T); \mU)} \int_0^T \left( \|Cx(t)\|_{\mH}^2 + \|u(t)\|_{\mU}^2  \right) dt,  
\ee
where $x(\cdot)$ is the solution of $x'(t) = A x(t) + Bu(t)$ on $(0, T)$ such that $x(0) = x_0 \in \mH$
($z=0$ here). Let $\tu_{T,\mathrm{opt}}$ be the optimal control solution of \eqref{optimal-control-finite}, and let $\tx_{T,\mathrm{opt}}$ be the corresponding solution and $\ty_{T,\mathrm{opt}}$ be the corresponding adjoint state. It is well known (see \cite[Section 9.2]{lasiecka2000control-v2}, see also \cite{FLT88}, \cite[Section 5]{Salamon_TAMS1987} and \cite{Flandoli87-LQR, LT91}) that there exists $P_T: [0, T] \to \cL(\mH)$ such that 
\be\label{state-adjointstate}
\ty_{T,\mathrm{opt}}(t) = P_T(t) \tx_{T,\mathrm{opt}}(t) \qquad\forall t \in [0, T]. 
\ee
When $B$ is bounded, \eqref{pro-technical-def-hg} and \eqref{rem-tech1} have been established in \cite{PZ13}, with $P$ replaced by $P_T$. 

Defining $h$ and $g$ by \eqref{pro-technical-def-hg} has an advantage on \eqref{state-adjointstate}. Indeed, thanks to \eqref{pro-technical-def-hg} where we use $P$ (instead of $P_T$), we can directly deal with $S_{\infty,\mathrm{opt}}^* (t)$ (see \eqref{pro-technical-pro-g}) instead of its approximation as in \cite{PZ13}. This also avoids to use the theory of differential and algebraic Riccati equations as in \cite{PZ13}.  

Last but not least, to derive the turnpike property from \Cref{pro-technical}, our analysis is different from the one performed in \cite{PZ13}. Instead of applying the Riccati theory, we apply a result (see \Cref{lem-transposition}) and follow some ideas of \cite{Ng-Riccati} to derive energy estimates (a kind of Lyapunov function), and then use a simple trick that is well known in the proof of the projection onto a closed convex set, using the parallelogram identity, to reach the conclusion.

\subsection{Organization of the paper} 
The paper is organized as follows. In \Cref{sect-LQ}, we revisit some facts of the linear quadratic optimal control theory in finite and infinite horizons. In \Cref{sect-stationary}, we establish the existence and uniqueness  of the optimal solutions and its Lagrange multipliers for the initial and approximate stationary optimization problems. The convergence of the approximate problem to the initial problem and  the proof of \Cref{pro-technical} are given in \Cref{sect-ingredient}.  \Cref{sect-main} is devoted to the proof of \Cref{thm-main}. The proof of several technical results are given in \Cref{sect-pro-Opt-T}. 

\medskip 
In what follows, for notational ease, we use $\langle \cdot, \cdot \rangle $ to denote $\langle \cdot, \cdot \rangle_{\mH}$ or $\langle \cdot, \cdot \rangle_{\mU}$ and $| \cdot |$  to denote $\| \cdot\|_{\mH}$ or $\| \cdot\|_{\mU}$ when the context is clear.

\section{Linear quadratic optimal control in finite and infinite horizons} \label{sect-LQ}
In this section, we first study a slightly more general setting than the one given by  \eqref{optimal-control-T} in finite horizon. We then derive its consequence for the corresponding linear quadratic optimal control problem in infinite horizon. 

Let $T>0$ and let $P_0 \in \cL(\mH)$ be nonnegative and symmetric. Given any $x \in C([0, T]; \mH)$ and $u \in L^2((0, T); \mU)$, we define the cost function
\be \label{def-JT}
J_T(x, u) = \int_0^T \left(|C x(t)|^2 + |u(t) |^2 \right) dt + \langle P_0 x(T), x(T) \rangle. 
\ee
Given any $\tau\in[0,T)$, any $x \in C([\tau, T]; \mH)$ and any $u \in L^2((\tau, T); \mU)$, we also define
$$
J_{\tau, T} (x, u) = \int_\tau^T \left(|C x(t) |^2 + |u(t) |^2 \right)  dt + \langle P_0 x(T), x(T) \rangle. 
$$
By strict convexity with respect to $u$, there exists a unique optimal control $\tu \in L^2((\tau, T); \mU)$ solution of 
$$
\inf_{u \in L^2((\tau, T); \mU)} J_{\tau, T} (x, u)
$$
where $x(\cdot)$ is the unique solution of $x'(t) = A x(t) + Bu(t)$ on $(\tau, T)$ such that $x(\tau) = \xi$. 
Since the problem is linear quadratic, $\tx$ and $\tu$ are linear functions of $\xi$. Moreover, there exists $c_{\tau, T} > 0$ such that the minimal value is bounded above by $c_T \| \xi \|_{\mH}^2$. Hence there exists a symmetric $P_T(\tau) \in \cL(\mH)$ such that  
\be  \label{def-Pt}
\langle P_T(\tau) \xi, \xi \rangle_{\mH} =  \inf_{u \in L^2((\tau, T); \mU)} J_{\tau, T} (x, u),  
\ee
Here, $P_T(T)(\xi, \xi)$ is understood as $\langle P_0 \xi, \xi \rangle $ for $\xi \in \mH$. Hence $\langle P_{T}(\tau)\xi, \xi \rangle$ is the cost on the time interval $[\tau,T]$ the initial data $\xi \in \mH$ at time $\tau$. 

Let us give more details on the concept of solution. Consider the slightly more general control system 
\be \label{CS-G}
\begin{split}
y'(t)&= Ay(t) + f(t) + Bu(t) + My(t) \quad  \mbox{on } (0, T), \\
y(0)&= y_0, 
\end{split}
\ee
with $y_0 \in \mH$, $f \in L^1((0, T); \mH)$, and $M \in \cL(\mH)$. A weak solution $y$ of \eqref{CS-G} is understood as an element $y \in C([0, T]; \mH)$ such that $y(0)=y_0$ and
\be\label{meaning-CS-G}
\frac{d}{dt} \langle y(t), \varphi \rangle_{\mH}  = \langle A y(t) + f(t) + Bu(t) + My(t), \varphi \rangle_{\mH} \quad\mbox{on } (0, T) \qquad\forall \varphi \in \cD(A^*) ,
\ee
for which the differential equation in \eqref{meaning-CS-G} is understood in the distributional sense, and the term $ \langle A y(t) + f(t) + Bu(t) + My(t), \varphi \rangle_{\mH} $ is understood as $\langle y(t), A^*\varphi \rangle_{\mH}  +  \langle f(t) + My(t), \varphi \rangle_{\mH} + \langle u(t), B^*\varphi \rangle_{\mU}$. 
The well-posedness of \eqref{CS-G} is known (see, e.g., \cite{BDDM07, CZ95, Ng-Riccati, WR00}). Concerning \eqref{CS-G}, we have the following result borrowed from \cite{Ng-Riccati}, which will be used repeatedly in this paper. 

\begin{lemma}\cite[Lemma 3.1]{Ng-Riccati}
\label{lem-transposition}
Let $T>0$, $y_0 \in \mH$,  $f \in L^1((0, T); \mH)$, $u \in L^2([0, T]; \mU)$, and $M \in \cL(\mH)$,  and let $y \in C([0, T]; \mH)$ be the unique weak solution of \eqref{CS-G}.  We have, for $t \in (0, T]$,  for $z_t \in \cD(A^*)$, and for $g \in C([0, t]; \cD(A^*))$,  
\begin{multline}\label{lem-transposition-cl1}
 \langle y(t), z_t \rangle_{\mH}  - \langle y_0, z(0) \rangle_{\mH} =  \int_0^t \langle u (s), B^*z(s) \rangle_{\mU} \, ds  \\ -  \int_0^t \langle y(s),  g(s) \rangle_{\mH} \, ds 
+  \int_0^t \langle f(s),  z(s) \rangle_{\mH} \, ds + \int_0^t \langle M y(s), z (s) \rangle_{\mH} \, ds, 
\end{multline}
where $z \in C([0, t]; \mH)$ is the unique weak solution of the backward system 
\be \label{lem-transposition-z}
\begin{split}
z'(s) &= - A^* z(s) - g(s)  \quad\mbox{on }  (0, t), \\
z(t) &= z_t.  
\end{split} 
\ee
Consequently, for $z_T \in \mH$ and $g \in L^1((0, T); \mH)$, the unique weak solution $z \in C([0, T]; \mH)$ of \eqref{lem-transposition-z} with $t = T$ satisfies
$$
\| B^*z \|_{L^2((0, T); \mU)} \le C_T \left(\| g\|_{L^1((0, T); \mH)} + \| z_T\|_{\mH} \right),
$$
and \eqref{lem-transposition-cl1} holds for $z_t \in \mH$ and $g \in L^1((0, t); \mH)$. Here $C_T>0$ does not depend on $g$, $f$, $z_T$. 
\end{lemma}

\begin{remark} \rm The equality \eqref{lem-transposition-cl1} can be seen as an integration by parts on time.   \Cref{lem-transposition} is related to solutions defined by transposition (see \cite[Remark 3.6]{Ng-Riccati} for related results). 
\end{remark}

The main result of this section is the following.

\begin{proposition} \label{pro-Opt-T} 
Let $T>0$ and let $P_0 \in \cL(\mH)$ be nonnegative and symmetric. Given $x_0, z \in \mH$, let $\tu$ be the optimal control solution of \eqref{optimal-control-T}. Let $\tx$ be the corresponding solution and $\ty$ be the corresponding adjoint state, i.e., 
\begin{equation} \label{pro-Opt-T-sys-opt1}
\begin{split}
& \tx'(t) = A \tx(t) + B \tu(t), \qquad \ty'(t) = - A^* \ty(t) - C^* (C \tx(t) - z)  \quad\mbox{ on } (0, T), \\
& \tx(0) = x_0, \qquad\qquad\qquad\quad \ty(T) = P_0 \tx(T) .
\end{split}
\end{equation} 
Then, we have
\be\label{pro-Opt-T-sys-opt1-u}
\tu = - B^* \ty  \quad\mbox{in } L^2((0, T); \mU).  
\ee 
If $z = 0$, then  
\be \label{pro-Opt-T-sys-opt2}
\ty (t) = P_T(t) \tx (t)  \qquad\forall t\in [0, T], 
\ee
where $P_T$ is defined by \eqref{def-Pt}.  As a consequence, when $B$ is bounded, we have
\be \label{pro-Opt-T-P}
P_{T}'(t) + A^*P_{T}(t) + P_{T}(t) A + C^* C - P_{T}(t)  B B^* P_{T}(t) = 0 \qquad\forall t\in [0, T], 
\ee
in the sense that 
\begin{multline*}
\frac{d}{dt} \langle P_T(t) \xi, \eta \rangle + \langle P_T(t) \xi, A \eta \rangle + \langle A \xi, P_T(t) \eta \rangle + \langle C \xi, C \eta \rangle \\
- \langle B^* P_T(t) \xi, B^* P_T(t) \eta \rangle = 0 \qquad\forall t\in [0, T],
\qquad \forall \xi, \eta \in \cD(A).
\end{multline*}
\end{proposition}

\begin{remark} \label{rem-B*-L2} \rm Since $\ty'(t) = - A^* \ty(t) - C^*(C \tx(t) - z) $ on $(0, T)$ and $\ty(T) \in \mH$, it follows from \Cref{lem-transposition} that 
$B^* \ty \in L^2((0, T); \mU) $ and 
\be \label{Ineq-Obs}
\| B^* \ty \|_{ L^2((0, T); \mU) } \le c_T \left( \|\ty(T)\|_{\mH} + \| C^*(C \tx - z) \|_{L^1((0, T); \mH)}  \right). 
\ee
Therefore, \eqref{pro-Opt-T-sys-opt1-u} makes sense. 
\end{remark}

\Cref{pro-Opt-T} is proved in \Cref{sect-pro-Opt-T}. 

\begin{remark} \label{rem-pro-Opt-T} \rm When $B$ is bounded and $z=0$, \Cref{pro-Opt-T} follows from Riccati theory (see, e.g.,  \cite{CZ95,BDDM07,Zabczyk08}). Assertions \eqref{pro-Opt-T-sys-opt1} and \eqref{pro-Opt-T-sys-opt1-u} are somehow known and are used to characterize the optimal solutions of \eqref{optimal-control-T}. Precisely, they are known when $z=0$, $P_0 = 0$, and $A$ enjoys some symmetric and coercivity properties (see \cite[Theorem 2.1 on page 114]{Lions71}), when $z=0$ and $A$ generates a strongly continuous {\it group} (see \cite[Theorem 2.3.1]{Flandoli87-LQR} in which the identity $\tu = - B^* \ty$ is mentioned to hold almost everywhere in $(0, T)$). A connection between $\tu$ and $\tx$ has been obtained in the setting considered here in \cite[Section 9.2]{lasiecka2000control-v2} (see also \cite[Theorem 2.1]{FLT88}, \cite[Section 5]{Salamon_TAMS1987} and \cite{Flandoli87-LQR,LT91}). In \cite{lasiecka2000control-v2}, $P_T(\tau)$ is first defined via a formula using $\tx$ and it is later proved that \eqref{def-Pt} holds. The proof that we give here is somehow in the same spirit but different from and more direct than the known ones mentioned above. Our proof also uses \Cref{lem-transposition}.
\end{remark}

\begin{remark} \rm The consideration of both states and adjoint states together are useful in the analysis of the paper. This was previously used 
to derive the rapid stabilisation of nonlinear control systems in several settings using Gramian operators (see \cite{Ng-Riccati,Ng-S-KdV,Ng-S-Schrodinger}). 
\end{remark}

We next present a useful consequence of \Cref{pro-Opt-T} for the optimal control problem in infinite horizon. 
It is obtained by applying \Cref{pro-Opt-T} with $z=0$ and by noting that, when $P_0=P$, we have $P_{T}(\tau) = P_0 = P$ for every $\tau \in [0, T]$.

\begin{corollary} \label{pro-P}  Assume that the triple $(A, B, C)$ enjoys the finite cost condition. Let $P \in \cL(\mH)$ be the symmetric operator defined by \eqref{lem-P-P} in \Cref{lem-P}. Let $T>0$. Given any $x_0 \in \mH$, let $\tu$ be the optimal control solution of
$$
\inf_{u \in L^2((0, T); \mU)} J_{T} (x, u)
$$
where $x(\cdot)$ is the unique solution of $x'(t) = A x(t) + Bu(t)$ on $(0, T)$ such that $x(0) = x_0$, and where  $J_T$ is defined by \eqref{def-JT} with $P_0 = P$. Let $\tx,\ty\in C((0, T); \mH)$ be the corresponding solution and adjoint state, i.e.,
\begin{equation*} 
\begin{split}
& \tx'(t) = A \tx(t) + B \tu(t), \qquad \ty'(t) = - A^* \ty(t) - C^*C \tx(t)  \quad\mbox{ on } (0, T), \\
& \tx(0) = x_0, \qquad\qquad\qquad\quad \ty(T) = P_0 \tx(T) .
\end{split}
\end{equation*} 
Then, we have 
$$
\tu = - B^* \ty  \quad\mbox{in } L^2((0, T); \mU), 
$$
and
\be \label{pro-P-sys-opt2}
\ty (t) = P \tx (t)  \qquad\forall t\in[0, T]. 
\ee
When $B$ is bounded, we have
\be \label{pro-P-P}
A^*P + P A + C^* C - P B B^* P = 0 
\ee
in the sense that
$$
\langle P \xi, A \eta \rangle + \langle A \xi, P \eta \rangle + \langle C \xi, C \eta \rangle - \langle B^* P \xi, B^* P \eta \rangle = 0 \qquad\forall \xi, \eta \in \cD(A). 
$$
\end{corollary}

\begin{remark} \rm A standard approach to derive the solution to an infinite-horizon optimal control problem is to consider the limit of the corresponding finite-horizon optimal problem with zero final cost. This is done in this way in \cite{FLT88,lasiecka2000control-v1} but the analysis is quite involved. We show here that the infinite-horizon optimal control problem can be more easily addressed by considering a finite-horizon problem with a suitably chosen final cost.
\end{remark}

\section{Stationary optimization problems}
\label{sect-stationary}

In this section, we establish \Cref{pro-stationary} in \Cref{sect-stationary-1}, state and prove a useful consequence of \Cref{pro-stationary} for the approximate problems in \Cref{sect-stationary-2} (the ones associated with $(A, B_k)$ for large positive $k$), and study the convergence to the initial problem in   \Cref{sect-stationary-3}. 

\subsection{Proof of \Cref{pro-stationary}} \label{sect-stationary-1} 
The existence of $(\bar x, \bar u)$ and the uniqueness of $(C\bar x, \bar u)$ are standard. Let us prove the uniqueness of $\bar x$. Assume that $(\bar x_1, \bar u), (\bar x_2, \bar u) \in \mH \times \mU$ are two optimal solutions. Setting $\xi = \bar x_2 - \bar x_1$, we have $A \xi = 0$ since $A \bar x_1 + B \bar u = A \bar x_2 + B \bar u=0$, and $C \xi = 0$ since $C \bar x_1 = C \bar x_2$.  We infer from \eqref{pro-stationary-AC} that $\xi =0$. This gives the uniqueness of $\bar x$. 

By the same arguments using the observability of $(A^*, B^*)$ and \eqref{pro-stationary-A*B*}, one can show that there exists at most one $\bar y \in \cD(A^*)$ such that \eqref{pro-stationary-cl2} holds. 

It thus remains to prove the existence of $\bar y$. We follow the ideas of \cite{PZ13} with some modifications and simplifications.  Let $u_{T,\mathrm{opt}}$ be the optimal control solution of \eqref{optimal-control-T} with $x_0 = \bar x$ and $P_0 = 0$. Let $x_{T,\mathrm{opt}}$ and $y_{T,\mathrm{opt}}$ be the corresponding solution and adjoint state. We have 
\begin{equation}\label{xopt0}
\begin{split}
x_{T,\mathrm{opt}}'(t) &= A x_{T,\mathrm{opt}}(t) + B u_{T,\mathrm{opt}}(t) \quad\mbox{on }(0, T),\\
x_{T,\mathrm{opt}}(0) &= \bar x ,
\end{split}
\end{equation}
and
\begin{equation}\label{yopt0}
\begin{split}
y_{T,\mathrm{opt}}'(t) &= - A^* y_{T,\mathrm{opt}}(t) - C^* (C x_{T,\mathrm{opt}}(t) - z) \quad\mbox{on } (0, T), \\
y_{T,\mathrm{opt}} (T) &= 0.     
\end{split}
\end{equation}
Since $A \bar x + B \bar u=0$, by optimality we must have
$$
\int_0^T \left(|C x_{T,\mathrm{opt}}(t) - z|^2 + |u_{T,\mathrm{opt}}(t)|^2 \right) dt \le T (|C \bar x - z|^2 + |\bar u|^2). 
$$
This implies, by Jensen's inequality, that
\be \label{pro-stationary-average-p1}
 \left|C \frac{1}{T}\int_0^T x_{T,\mathrm{opt}}(t) \, dt  - z\right|^2 + \left|\frac{1}{T}\int_0^T u_{T,\mathrm{opt}}(t) \, dt\right|^2  \le |C \bar x - z|^2 + |\bar u|^2. 
\ee
Using the arguments in \cite[Remark 2.1]{PZ13}, the finite-time observability of $(A, C)$ and \eqref{pro-stationary-average-p1} imply that
\be \label{pro-stationary-coucou}
\int_0^T |x_{T,\mathrm{opt}}(t)|^2\, dt \le c T. 
\ee
Here and in what follows in this proof, $c$ denotes a generic positive constant independent of $T$ and of $s \in [0, T]$.
It follows from \eqref{pro-stationary-average-p1} that there exists a sequence $T_n \to + \infty$ such that 
$$
\frac{1}{T_n}\int_0^{T_n} x_{T_n, opt} \rightharpoonup \tx \mbox{ weakly in } \mH \quad \mbox{ and } \quad  \frac{1}{T_n}\int_0^{T_n} u_{T_n, opt} \rightharpoonup \tu \mbox{ weakly in } \mU,  
$$
for some $\tx \in \mH$ and $\tu \in \mU$ \footnote{The notation $\rightharpoonup$ means the weak convergence.}. We derive from \eqref{pro-stationary-average-p1} and the properties of the weak convergence that 
\be \label{pro-stationary-txtu-p1}
|C \tx - z|^2 + |\tu|^2 \le |C \bar x - z|^2 + |\bar u|^2. 
\ee
Integrating \eqref{xopt0}, we obtain 
\be\label{pro-stationary-p1}
\frac{x_{T,\mathrm{opt}}(T) - x_{T,\mathrm{opt}}(0)}{T} = A \frac{1}{T}\int_0^T x_{T,\mathrm{opt}}  + B \frac{1}{T}\int_0^T u_{T,\mathrm{opt}} . 
\ee
Since, by  \eqref{pro-stationary-coucou},  
\be\label{pro-stationary-AC1}
| x_{T,\mathrm{opt}}(T)| \le c T^{1/2},
\ee
 we then infer from \eqref{pro-stationary-p1} that 
\be \label{pro-stationary-txtu-p2}
A \tx + B \tu = 0. 
\ee
Combining \eqref{pro-stationary-txtu-p1} and \eqref{pro-stationary-txtu-p2} and the uniqueness of $(\bar x, \bar u)$, we infer that $\tx = \bar x$ and $\tu = \bar u$. Moreover, 
\be\label{pro-stationary-hihi}
C \frac{1}{T}\int_0^T x_{T,\mathrm{opt}}   \to C \bar x \quad \mbox{ and } \quad  \frac{1}{T}\int_0^T u_{T,\mathrm{opt}} \to \bar u
\ee
as $T \to + \infty$, with a strong convergence. 

Integrating \eqref{yopt0}, we have 
\be \label{pro-stationary-y-p1}
\frac{y_{T,\mathrm{opt}}(T) - y_{T,\mathrm{opt}}(0)}{T} = - A^* \frac{1}{T}\int_0^T y_{T,\mathrm{opt}}   - C^* \left(C \frac{1}{T}\int_0^T x_{T,\mathrm{opt}}  - z \right). 
\ee
Since $(A^*, B^*)$ is finite-time observable, as in the proof of \eqref{pro-stationary-AC1}, we infer that $|y_{T,\mathrm{opt}}(0)| \le c T^{1/2}$. We deduce that the left-hand side of \eqref{pro-stationary-y-p1} converges to $0$ as $T \to + \infty$. Without loss of generality, we assume that  
$$
 \frac{1}{T_n}\int_0^{T_n} y_{T_n, opt}   \rightharpoonup  \bar y \quad\mbox{ in } \mH,  
$$
for some $\bar y \in \mH$.  Using  \eqref{pro-stationary-hihi}, we then infer from \eqref{pro-stationary-y-p1} that $\bar y \in \cD(A^*)$, that $-A^* \bar y  =  C^* (C \bar x - z)$, and that
$$
 \frac{1}{T_n}\int_0^{T_n} y_{T_n, opt}   \rightharpoonup  \bar y \quad\mbox{weakly  in } \cD(A^*). 
$$
Since $u_{T,\mathrm{opt}}=- B^* y_{T,\mathrm{opt}}$ by \Cref{pro-Opt-T} which yields 
$ -B^* \frac{1}{T}\int_0^T y_{T,\mathrm{opt}} = \frac{1}{T}\int_0^T u_{T,\mathrm{opt}}$, we finally obtain 
$-B^* \bar y =  \bar u$. The existence of $\bar y$ is proved. The proof is complete.

\begin{remark} \rm The above proof is in the spirit of the one given in \cite{PZ13}. Nevertheless, instead of considering an arbitrary initial data in the definition of $u_{T,\mathrm{opt}}$ as in \cite{PZ13}, we have considered the initial data $\bar x$, noting that $A \bar x + B \bar u=0$, thus allowing us to compare the cost of $(x_{T,\mathrm{opt}}, u_{T,\mathrm{opt}})$ with the cost of $(\bar x, \bar u)$. As a result, the arguments are simpler and require weaker assumptions. 
\end{remark}

\subsection{Consequence of \Cref{pro-stationary} for the approximate stationary problems} \label{sect-stationary-2}

A useful consequence of \Cref{pro-stationary} for the approximate stationary problems associated with $(A, B_k)$ is now stated. 

\begin{corollary}\label{cor-stationary} Let $z \in \mH$.  Assume that $(A, C)$ is finite-time observable in time $T_0>0$ and that $\ker A \cap \ker C = \{0 \}$. For any $k\in\N$ large enough, recalling that $B_k$ is defined by \eqref{def-Bk}, there exists a unique $(\bar u_k, \bar x_k) \in \mU \times \mH$ such that $A \bar x_k + B_k \bar u_k = 0$ and 
$$
\|C \bar x_k - z\|_{\mH}^2 + \|\bar u_k\|_{\mU}^2 = \inf_{(x,u)\in\mH\times\mU\atop Ax + B_k u = 0} \left( \|Cx - z\|_{\mH}^2 + \|u\|_{\mU}^2 \right) . 
$$
Assume in addition that $(A^*, B^*)$ is finite-time observable in time $T_0>0$ and that $\ker A^* \cap \ker B^*= \{0 \}$. Then there exists a unique $\bar y_k \in \cD(A^*)$ such that 
\be \label{cor-stationary-cl2}
- A^* \bar y_k - C^* (C \bar x_k - z) = 0 \quad \mbox{ and } \quad \bar u_k = - B_k^* \bar y_k. 
\ee
\end{corollary}

\begin{proof} We first observe that if $\ker A^* \cap \ker B^* = \{ 0\}$ 
then $\ker A^* \cap \ker B_k^* = \{0 \}$.
Indeed, let $\xi \in \ker A^* \cap \ker B_k^*$. Since $
A^* \xi = 0$, it follows that 
$A^* J_k^* \xi = J_k^* A^* \xi = 0$. 
On the other hand,  $0 = B_k^* \xi = B^* J_k^* \xi$.
Since $\ker A^* \cap \ker B^* = \{ 0\}$, we infer that $J_k^* \xi = 0$ and therefore $\xi =0$. The observation is proved. 

As a consequence of \Cref{pro-stationary}, it now suffices to note that $(A, B_k)$ is finite-time observable in time $T_0$ if $(A, B)$ is finite-time observable in time $T_0$. This is so because, for $\xi \in \mH$,  
$$
\int_0^{T_0} \|B_k^* e^{t A^*} \xi  \|^2 \, dt = \int_0^{T_0} \|B^* J_k^* e^{t A^*} \xi  \|^2 \, dt =  \int_0^{T_0} \|B^* e^{t A^*} J_k^* \xi  \|^2 \, dt  \ge c \| J_k^* \xi \|_{\mH}^2 \ge c_k \| \xi \|_{\mH}, 
$$
for some positive constant $c_k$ independent of $\xi$. 
\end{proof}

\subsection{Convergence properties for the approximate stationary problems}\label{sect-stationary-3}

\begin{lemma} \label{lem1} 
Under Assumption \ref{H2}, as $k \to + \infty$, 
$$
\bar x_k \to \bar x \mbox{ in } \mH, \quad \bar u_k \to \bar u \mbox{ in } \mU, \quad \mbox{ and } \quad 
\bar y_k \to \bar y \mbox{ in } \mH. 
$$
Here $(\bar x_k, \bar u_k, \bar y_k)$ and $(\bar x, \bar u, \bar y)$ are defined in \Cref{cor-stationary} and \Cref{pro-stationary}, respectively. 
\end{lemma}

\begin{proof}  By optimality of $(\bar x_k, \bar u_k)$, we have 
\be \label{lem1-PE}
|C \bar x_k - z|^2 + |\bar u_k|^2 \le |z|^2. 
\ee
Noting that $A \bar x_k + B_k \bar u_k=0$, since $(A, C)$ is finite-time observable, we infer from \eqref{lem1-PE} that
$|\bar x_k|^2  \le c |z|^2$.
Let $(\bar x_{n_k})$ and $\bar u_{n_k}$ be a subsequence of $(\bar x_k)$ and $(\bar u_k)$ such that, for some $(\hx, \hu) \in \mH \times \mU$,
\be \label{lem1-xkuk}
\bar x_{n_k} \rightharpoonup \hx \mbox{ weakly in $\mH$} \quad \mbox{ and } \quad  \bar u_{n_k} \rightharpoonup \hu \mbox{ weakly in $\mU$}.
\ee
Then 
\be \label{lem1-xkuk-p1}
A \hx + B \hu = 0, \quad \mbox{ and } \quad \liminf_{k \to + \infty} \big( |C \bar x_{n_k} - z|^2 + |\bar u_{n_k}|^2 \big) \ge 
 |C  \hx - z|^2 + |\hu|^2. 
\ee
On the other hand, since $A \hx +  B  \hu =0$, we have 
$A (J_k \hx) + B_k \hu = J_k(A \hx + B \hu) = 0$.
This implies, by optimality of $(\bar x_k, \bar u_k)$, 
\be \label{lem1-p2-1}
 |C J_k \hx - z|^2 + |\hu|^2 \ge  |C \bar x_k - z|^2 + |\bar u_k|^2. 
\ee
Since $ |C J_k \hx - z|^2 + |\hu|^2  \to  |C \hx - z|^2 + |\hu|^2$ as $k \to + \infty$, we infer from \eqref{lem1-p2-1} that 
\be\label{lem1-p2}
\limsup_{k \to + \infty} \big( |C \bar x_k - z|^2 + |\bar u_k|^2 \big) \le |C \hx - z|^2 + |\hu|^2. 
\ee
As in the proof of \eqref{lem1-p2} we also have, since $A \bar x + B \bar u = 0$,  
\be\label{lem1-p2-*}
\limsup_{k \to + \infty} \big( |C \bar x_k - z|^2 + |\bar u_k|^2 \big) \le |C \bar x - z|^2 + |\bar u|^2. 
\ee
Combining \eqref{lem1-xkuk-p1} and \eqref{lem1-p2} yields 
\be\label{lem1-xkuk-p3}
\lim_{k \to + \infty} \big( |C \bar x_{n_k} - z|^2 + |\bar u_{n_k}|^2 \big) = |C \hx - z|^2 + |\hu|^2. 
\ee
From \eqref{lem1-p2-*} and \eqref{lem1-xkuk-p3}, we infer that $|C \hx - z|^2 + |\hu|^2 \le |C \bar x - z|^2 + |\bar u|^2$, and therefore, by optimality of $(\bar x, \bar u)$, the equality holds.
By uniqueness of $(\bar x, \bar u)$, we thus have
\be \label{lem1-xkuk-p4-*}
(\hx, \hu) = (\bar x, \bar u). 
\ee
From \eqref{lem1-xkuk}, \eqref{lem1-xkuk-p3}, and \eqref{lem1-xkuk-p4-*}, we infer that 
$$
\bar u_{k} \to \bar u \mbox{ in } \mU \quad \mbox{ and } \quad C \bar x_{k} \to C \bar x \mbox{ in } \mH
$$
with strong convergence and for the whole sequence.

We next prove that 
$$
\bar x_k \to \bar x \mbox{ in } \mH. 
$$
Since $A \bar x_k + B_k \bar u_k = 0$, we have
\be \label{lem1-cc1}
\bar x_k = e^{t A} \bar x_k + J_k f_k (t) \mbox{ for } t \ge 0 \quad \mbox{ where } f_k (t) = \int_0^t e^{(t-s)A} B \bar u_k \, ds ,
\ee
and since $A \bar x  + B \bar u = 0$, we obtain 
\be \label{lem1-cc2}
\bar x = e^{t A} \bar x  + f(t) \mbox{ for } t \ge 0 \quad\mbox{ where } f (t) = \int_0^t e^{(t-s)A} B \bar u \, ds .
\ee
Since $u_k \to u $ in $\mH$, we infer from the admissibility of $B$ that, for every $T>0$, 
\be  \label{lem1-cc3}
f_k \to f \quad\mbox{ in } C([0, T]; \mH), 
\ee
which shows that
\be \label{lem1-cc4}
C J_k f_k \to C f \quad\mbox{ in } C([0, T]; \mH).  
\ee
Combining \eqref{lem1-cc1}, \eqref{lem1-cc2}, and \eqref{lem1-cc4}, and using the fact that $C \bar x_k \to C \bar x  $ in $\mH$, we derive that 
\be  \label{lem1-cc5} 
\left(t\mapsto C e^{t A} \bar x_k - C e^{t A} \bar x\right)\to 0 \quad\mbox{ in } C ([0, T]; \mH). 
\ee
Since $(A, C)$ is finite-time observable in time $T_0$, we infer from \eqref{lem1-cc5} that
\be  \label{lem1-cc6}
e^{T_0 A} (\bar x_k - \hat x) \to 0 \quad\mbox{ in } \mH.  
\ee
Combining \eqref{lem1-cc1} and \eqref{lem1-cc2} with $t = T_0$, and using \eqref{lem1-cc3} and \eqref{lem1-cc6}, we get 
\be \label{lem1-xkuk-cl}
\bar x_k \to \bar x \mbox{ in } \mH. 
\ee
It thus remains to show that $\bar y_k \to \bar y \mbox{ in } \mH$.
Since
$$
A^* \bar y \mathop{=}^{\eqref{pro-stationary-cl2}} - C^* (C \bar x - z), \quad A^* \bar y_k \mathop{=}^{\eqref{cor-stationary-cl2}} - C^* (C \bar x_k - z), \quad \mbox{ and } \quad C \bar x_k \mathop{\to}^{\eqref{lem1-xkuk-cl}} C \bar x \mbox{ in } \mH, 
$$
and
$$
B^* \bar y \mathop{=}^{\eqref{pro-stationary-cl2}} \bar u, \quad B^* J_k^* \bar y_k \mathop{=}^{\eqref{cor-stationary-cl2}} \bar u_k, \quad \mbox{ and } \quad u_k \mathop{\to}^{\eqref{lem1-xkuk-cl}} u \mbox{ in } \mU,   
$$
we infer that, as $k \to + \infty$, 
$$
A^* J_k^* \bar y_k = J_k^* A^*  \bar y_k \to A^* \bar y \mbox{ in } \mH \quad \mbox{ and } \quad B^* J_k^*y_k \to B^* y \mbox{ in } \mU. 
$$
By the finite-time observability of $(A^*, B^*)$, using the arguments in \cite[Remark 2.1]{PZ13}, we deduce that $J_k^* \bar y_k \to \bar y \mbox{ in } \mH$ as $k \to + \infty$.
Since  $A^* J_k^* \bar y_k \to A^* \bar y$ in $\mH$, it follows that 
$J_k^* \bar y_k -  \bar y \to 0 \mbox{ in } \cD(A^*)$ as $k \to + \infty$.
This implies that $(J_k^*)^{-1} (J_k^* \bar y_k -  \bar y)  \to 0$ in $\mH$.
Since $\bar y \in \cD(A^*)$,  we finally obtain that $\bar y_k \to \bar y \mbox{ in } \mH$ as $k \to + \infty$.
The proof is complete. 
\end{proof}

\section{Proof of \Cref{pro-technical}}
\label{sect-ingredient}

This section is devoted to the proof of \Cref{pro-technical} and consists of two subsections. In the first one, we establish several results used in the proof of \Cref{pro-technical}. The proof of \Cref{pro-technical} is given in the second subsection.

\subsection{Some useful lemmas}

We first establish \Cref{pro-technical} when $B$ is bounded. 

\begin{lemma}\label{lem-technical1} Under Assumption \ref{H2} of \Cref{thm-main}, let $z \in \mH$ and $T>0$. Given $x_0 \in \mH$, define $h$ by
\be \label{lem-technical1-def-hg}
h(t) = y_{T,\mathrm{opt}}(t) - \bar y  -  P (x_{T,\mathrm{opt}}(t) - \bar x) \qquad\forall t\in [0, T], 
\ee
where $P$ is defined by \eqref{lem-P-P}, $x_{T,\mathrm{opt}}$ is the corresponding optimal solution of \eqref{optimal-control-T},  $y_{T,\mathrm{opt}}$ is the corresponding adjoint state, and $(\bar x, \bar y)$ is given in \Cref{pro-stationary}.  Assume in addition that $B$ is bounded. Then 
\be\label{lem-technical1-h}
h'(t) = \left(- A^* + P B B^* \right) h(t) \qquad\forall t\in [0, T].
\ee
Consequently, setting $g(t) = h(T - t)$ for every $t\in[0, T]$, we have
$$
g(t) = S_{\infty,\mathrm{opt}}^* (t) g(0) \qquad\forall t\in [0, T],
$$
where $S_{\infty,\mathrm{opt}}(t)$ is defined in \Cref{pro-def-Sopt} and $S_{\infty,\mathrm{opt}}^*(t)$ is its adjoint.
\end{lemma}

\begin{proof} 
We infer from \Cref{pro-Opt-T} that, on $(0,T)$,
\be \label{lem-technical-Sys}
\begin{split}
x_{T,\mathrm{opt}}'(t)&= A x_{T,\mathrm{opt}}(t) - B B^* y_{T,\mathrm{opt}}(t) ,\qquad\qquad\qquad x_{T,\mathrm{opt}}(0) = x_0, \\
y_{T,\mathrm{opt}}'(t) &= - A^* y_{T,\mathrm{opt}}(t) - C^* (C x_{T,\mathrm{opt}}(t) - z) , \qquad\ y_{T,\mathrm{opt}} (T) = 0. 
\end{split}
\ee
Using \eqref{pro-stationary-cl1} and \eqref{pro-stationary-cl2}, we derive from \eqref{lem-technical-Sys} that, on $(0,T)$, 
\be \label{lem-technical-Sys-mod}
\begin{split}
(x_{T,\mathrm{opt}} - \bar x)'(t) &= A (x_{T,\mathrm{opt}}(t) - \bar x)  - B B^* (y_{T,\mathrm{opt}}(t) - \bar y), \\
(y_{T,\mathrm{opt}} - \bar y)'(t) &= - A^* (y_{T,\mathrm{opt}}(t) - \bar y)  - C^* C (x_{T,\mathrm{opt}}(t) -  \bar x) .
\end{split}
\ee
Since $B$ is bounded, $(\bar x, \bar u) \in \mH \times \mU$, and $A \bar x + B \bar u = 0$, we infer that $\bar x \in \cD(A)$. By definition of $h$ in \eqref{lem-technical1-def-hg}, we formally have, on $(0, T)$,
\begin{equation}\label{lem-technical1-h1}
\begin{split}
h'(t) &= (y_{T,\mathrm{opt}}(t) - \bar y)' - P (x_{T,\mathrm{opt}}(t) - \bar x) '  \\
&\mathop{=}^{\eqref{lem-technical-Sys-mod}}  - A^* (y_{T,\mathrm{opt}}(t) - \bar y) - C^* C (x_{T,\mathrm{opt}}(t) - \bar x) - P\left( A (x_{T,\mathrm{opt}}(t) - \bar x) - B B^* (y_{T,\mathrm{opt}}(t) - \bar y) \right) .
\end{split}
\end{equation}
Using \eqref{pro-P-P}, we infer from \eqref{lem-technical1-h1} that, on $(0, T)$,
$$
h'(t) = (-A^* + P B B^*) (y_{T,\mathrm{opt}}(t) - \bar y) - (-A^* P + P B B^* P ) (x_{T,\mathrm{opt}}(t) - \bar x) ,
$$
which yields \eqref{lem-technical1-h} by the definition of $h$ in \eqref{lem-technical1-def-hg}. 

The rigorous proof of \eqref{lem-technical1-h} goes as follows. We first assume that $x_0 \in \cD(A)$. Let $\xi_n \in C([0, T]; \cD(A))$ and $\eta_n \in C([0, T]; \cD(A))$ be such that 
$$
\xi_n \to C^* C(x_{T,\mathrm{opt}} - \bar x) \mbox{ in } C([0, T]; \mH) \quad \mbox{ and } \quad \eta_n \to  B B^*(y_{T,\mathrm{opt}} - \bar y) \mbox{ in } C([0, T]; \mH). 
$$
Let $x_n$ and $y_n$ in $C([0, T]; \mH)$ be the solutions of the systems
\begin{equation*}
\left\{\begin{array}{c}
x_n' = A x_n  - \eta_n \mbox{ in } (0, T), \\
x_n(0)=x_0 - \bar x,
\end{array} \right. 
\quad \mbox{ and } \quad 
\left\{\begin{array}{c}
y_n' = - A^* y_n  -\xi_n \mbox{ in } (0, T),\\
y_n(T)=0, 
\end{array} \right. 
\end{equation*}
and set 
$h_n = y_n -  P x_n$ on $[0, T]$.
Then $x_n, y_n \in C^1([0,T]; \mH) \cap C([0, T]; \cD(A))$, and $x_n \to x_{T,\mathrm{opt}}- \bar x$ and $y_n \to y_{T,\mathrm{opt}}- \bar y$ in 
$C([0, T]; \mH)$. 
We have 
\be\label{lem-technical1-h1-11}
h_n' = y_n' - P x_n' 
= - A^* y_n - \xi_n -
P\big( A x_n - \eta_n \big) \mbox{ in } (0,T).
\ee
Using \eqref{pro-P-P}, we infer from \eqref{lem-technical1-h1-11} that, on $(0,T)$,
\begin{equation}\label{lem-technical1-hn}
\begin{split}
h_n' &= (-A^* + P B B^*) y_n  - (-A^* P + P B B^* P ) x_n -  P (B B^* y_n - \eta_n) + (C^* C x_n - \xi_n) \\
&=
 (-A^* + P B B^*) h_n - P (B B^* y_n - \eta_n) + (C^* C x_n - \xi_n) 
\end{split}
\end{equation}
Note that, as $n \to + \infty$. 
$$
-P (B B^* y_n - \eta_n) + (C^* C x_n - \xi_n)   \to 0 \mbox{ in } C([0,T];\mH) \quad \mbox{ and }\quad h_n \to h \mbox{ in } C([0,T];\mH).
$$
Letting $n \to + \infty$ in \eqref{lem-technical1-hn}, we obtain \eqref{lem-technical1-h}. The proof of \eqref{lem-technical1-h} in the case $x_0 \in \mH$, follows from the case $x_0 \in \cD(A)$ by a standard approximation argument. 

It follows from \eqref{lem-technical1-h} that 
$g' = (A^* - P B B^*) g$ on $(0, T)$.
Noting that 
$(A^* - P B B^*)^* = A - B B^* P$, we finally infer that 
$g(t) = S_{\infty,\mathrm{opt}}^* (t) g(0)$ for every $t\in[0,T]$.
The proof is complete. 
\end{proof}

In the same spirit, we have the following useful result. 

\begin{lemma}\label{lem-technical2} Under Assumption \ref{H2} of \Cref{thm-main}, let $z \in \mH$ and $T>0$, and let $P_0 \in \cL(\mH)$ be symmetric and nonnegative. 
Given $x_0 \in \mH$, set 
\be \label{lem-technical2-def-hg}
h_T(t) = y_{T,\mathrm{opt}}(t) - \bar y  -  P_T(t) (x_{T,\mathrm{opt}}(t) - \bar x) \qquad\forall t\in [0, T],  
\ee
where $P_T(t)$ is defined by \eqref{def-Pt}, $x_{T,\mathrm{opt}}$ is the corresponding optimal solution of \eqref{optimal-control-T},  $y_{T,\mathrm{opt}}$ is the corresponding state, and $(\bar x, \bar y)$ is given in \Cref{pro-stationary}. Assume that $B$ is bounded. Then 
$$
h_T'  = \big(-A^* + P_T B B^* \big) g_T \quad\mbox{ on } (0, T). 
$$
\end{lemma}

\begin{proof}
By definition of $h_T$ in \eqref{lem-technical2-def-hg}, we formally have, in $(0, T)$,  
\begin{equation}\label{thm-main-eq-h1}
\begin{split}
h_T'(t) &= (y_{T,\mathrm{opt}}(t) - \bar y)' - P_T (x_{T,\mathrm{opt}}(t) - \bar x) '  - P_T' (x_{T,\mathrm{opt}}(t) - \bar x)  \\
&\mathop{=}^{\eqref{lem-technical-Sys-mod}}  - A^* (y_{T,\mathrm{opt}}(t) - \bar y) - C^* C (x_{T,\mathrm{opt}}(t) - \bar x) \\
& \qquad\qquad\qquad - P_T \big( A (x_{T,\mathrm{opt}}(t) - \bar x) - B B^* (y_{T,\mathrm{opt}}(t) - \bar y) \Big) - P_T' (x_{T,\mathrm{opt}}(t) - \bar x).
\end{split}
\end{equation}
Using \eqref{pro-Opt-T-P} of \Cref{pro-Opt-T}, we infer from \eqref{thm-main-eq-h1} that, on $(0,T)$,
$$
h_T'(t) = (-A^* + P_T B B^*) (y_{T,\mathrm{opt}}(t) - \bar y) - (-A^* P_T + P_T B B^* P_T ) (x_{T,\mathrm{opt}}(t) - \bar x),
$$
which yields, by the definition of $h_T$ in \eqref{lem-technical2-def-hg}, that
$h_T' = (-A^* + P_T B B^*) h_T$ on $(0, T)$.
The rigorous proof is done with an approximation argument similar to the one used to establish \eqref{lem-technical1-h}. The details are omitted. 
\end{proof}

We next study the convergence of the approximate problems to the initial one. The following result is the dynamical version of \Cref{lem1}, which was for the stationary problem.  

\begin{lemma} \label{lem-ap} Let $P_0 \in \cL(\mH)$ be symmetric and nonnegative, let $z \in \mH$ and $T > 0$. For large positive $k$, define $B_k$ by \eqref{def-Bk}. We have, for $\xi \in \mH$,  
\begin{multline} \label{lem-ap-cl1}
\lim_{k \to + \infty} \inf_{u \in L^2((0, T); \mU)} \left( \int_0^T \left(|C x_k(t) - z|^2 + |u(t) |^2 \right) dt + \langle P_0 x_k(T), x_k(T) \rangle \right) \\
= \inf_{u \in L^2((0, T); \mU)} \left( \int_0^T \left( |C x(t) - z|^2 + |u(t) |^2 \right) dt + \langle P_0 x(T), x(T) \rangle \right), 
\end{multline}
where 
\begin{align}
& x_k'(t) = A x_k(t) + B_k u(t) \quad\mbox{on }(0, T), \qquad x_k(0) = \xi , \label{lem-ap-xk} \\
& x'(t) = A x(t) + B u(t)\qquad\mbox{on }(0, T), \qquad x(0) = \xi. \label{lem-ap-x}
\end{align}
Let $u_{T,\mathrm{opt}, k}$ and $u_{T,\mathrm{opt}}$ be the optimal controls,  
let $x_{T,\mathrm{opt}, k}$ and $x_{T,\mathrm{opt}}$ be the corresponding solutions, and let
$y_{T,\mathrm{opt}, k}$ and $y_{T,\mathrm{opt}}$ be the corresponding adjoint states, respectively. 
We have 
\begin{align}
&u_{T,\mathrm{opt}, k} \to   u_{T,\mathrm{opt}} \quad\mbox{ in } L^2((0, T); \mU),  \label{lem-ap-cl2}\\
&x_{T,\mathrm{opt}, k} \to   x_{T,\mathrm{opt}} \quad\mbox{ in } C([0, T]; \mH), \label{lem-ap-cl3} \\
&y_{T,\mathrm{opt}, k} \to   y_{T,\mathrm{opt}} \quad\mbox{ in } C([0, T]; \mH). \label{lem-ap-cl4}
\end{align}

\end{lemma}

\begin{proof} Let $\tx_k \in C([0, T]; \mH)$ be the unique solution of 
\be \label{lem-ap-p1-0}
\tx_k'(t) = A \tx_k(t) + B_k u_{T,\mathrm{opt}}(t) \quad \mbox{ on } (0, T) \quad \mbox{ and } \quad \tx_k(0) = \xi
\ee
(note that we use here $u_{T,\mathrm{opt}}$ and not $u_{T,\mathrm{opt}, k}$). Then, for $t \in [0, T]$,  
$$
\tx_k(t) = e^{t A} \xi + \int_0^t e^{(t-s)A} B_k u_{T,\mathrm{opt}}(s) \, ds = e^{t A} \xi + J_k \int_0^t e^{(t-s)A} B u_{T,\mathrm{opt}}(s) \, ds. 
$$
Since $\big\{\int_0^t e^{(t-s)A} B u_{T,\mathrm{opt}}(s) \, ds \in \mH\ \mid\ t \in [0, T] \big\}$ is a compact subset of $\mH$ (by admissibility of $B$), we infer that $\tx_k \to x_{T,\mathrm{opt}}$ in $C([0, T]; \mH)$ as $k \to + \infty$.  This implies that
\begin{multline}\label{lem-ap-txk1}
\limsup_{k \to + \infty} \left( \int_0^T \left(|C \tx_k(t)- z|^2 + |u_{T,\mathrm{opt}}(t) |^2 \right) dt + \langle P_0 \tx_k(T), \tx_k(T) \rangle \right) \\
=  \int_0^T \left(|C x_{T,\mathrm{opt}}(t) - z|^2 + |u_{T,\mathrm{opt}}(t) |^2 \right) dt + \langle P_0 x_{T,\mathrm{opt}}(T), x_{T,\mathrm{opt}}(T) \rangle. 
\end{multline}
We deduce that 
\begin{multline} \label{lem-ap-p1-2}
\limsup_{k \to + \infty} \inf_{u \in L^2((0, T); \mU)} \int_0^T \left( |C x_k(t) - z |^2 + |u(t) |^2 \right) \, dt + \langle P_0 x_k(T), x_k(T) \rangle \\
\le \inf_{u \in L^2((0, T); \mU)} \int_0^T \left( |C x(t) - z |^2 + |u(t) |^2 \right) \, dt + \langle P_0 x(T), x(T) \rangle .
\end{multline}
where $x_k$ and $x$ are defined by \eqref{lem-ap-xk} and \eqref{lem-ap-x}. 

On the other hand, from \eqref{lem-ap-p1-2} we have
$$
\int_0^T \left( |C x_{T,\mathrm{opt}, k}(t)|^2 +  |u_{T,\mathrm{opt}, k}(t)|^2 \right) dt   \le c. 
$$
Using the finite-time observability property of $(A, C)$, we infer that 
$$
\int_0^T \left( |x_{T,\mathrm{opt}, k}(t)|^2 +  |u_{T,\mathrm{opt}, k}(t)|^2 \right) dt  \le c. 
$$
Without loss of generality, we assume that 
\be
x_{T,\mathrm{opt}, k} \rightharpoonup \tx \ \mbox{ weakly in } L^2((0, T);\mH) \quad \mbox{ and } \quad
u_{T,\mathrm{opt}, k} \rightharpoonup \tu \ \mbox{ weakly in } L^2((0, T);\mU),  
\ee
for some $\tx \mbox{ in } L^2((0, T);\mH)$ and $\tu \mbox{ in } L^2((0, T);\mU)$.  We infer that 
$$
\tx(t) = e^{t A} \xi + \int_0^t e^{(t-s)A} B \tu(s) \, ds  \qquad\forall t\in [0, T]. 
$$
This implies that
\begin{multline} \label{lem-ap-p1-3}
\inf_{u \in L^2((0, T); \mU)} \int_0^T \left(|C x(t) - z|^2 + |u(t) |^2 \right) \, dt + \langle P_0 x(T), x(T) \rangle \\
\le \liminf_{k \to + \infty} \inf_{u \in L^2((0, T); \mU)} \int_0^T \left(|C x_k(t) - z |^2 + |u(t) |^2 \right) \, dt + \langle P_0 x_k(T), x_k(T) \rangle, 
\end{multline}
where $x_k$ and $x$ are given by \eqref{lem-ap-xk} and \eqref{lem-ap-x}. 

Assertion \eqref{lem-ap-cl1} now follows from \eqref{lem-ap-p1-2} and \eqref{lem-ap-p1-3}.

\medskip
We next prove \eqref{lem-ap-cl2}, \eqref{lem-ap-cl3} and \eqref{lem-ap-cl4}. 
Given a Hilbert space $\mX$ and $a: \mX \times \mX \to \mC$ a bilinear Hermitian form, we recall the parallelogram identity
\be \label{thm-main-Parallelogram}
a\big( (f- g)/2, (f- g)/2\big) + a\big(  (f +  g)/2, (f +  g)/2\big)  = \frac{1}{2} \big(a(f, f) + a(g, g)\big) \qquad\forall f, g \in \mX.  
\ee
Using \eqref{thm-main-Parallelogram} with $\mX = L^2((0, T); \mH \times \mU)$ and  
$$
a (f_1, f_2) = \int_0^T \left( \langle C x_1(t), C x_2(t) \rangle_{\mH}  + \langle u_1(t), u_2(t) \rangle_{\mU} \right)  dt + \langle P_0 x_1(T), x_2(T) \rangle 
$$
with $f_1 = (x_1, u_1)$ and $f_2 = (x_2, u_2)$, 
as in the standard proof of the projection onto a closed convex subset of a Hilbert space, we infer from \eqref{lem-ap-cl1} and \eqref{lem-ap-txk1} (recall that  $\tx_k$ is defined by \eqref{lem-ap-p1-0}) that
\begin{multline*}
 \limsup_{k \to + \infty} \bigg( \int_0^T \left(|C (\tx_k(t)- x_{T,\mathrm{opt}, k}(t)) |^2 + |u_{T,\mathrm{opt}}(t) - u_{T,\mathrm{opt}, k}(t) |^2 \right)  dt  \\
+ \langle P_0 (\tx_k - x_{T,\mathrm{opt}, k})(T), (\tx_k - x_{T,\mathrm{opt}, k})(T) \rangle \bigg)  = 0. 
\end{multline*}
Therefore
$$
u_{T,\mathrm{opt}, k} \to u_{T,\mathrm{opt}} \mbox{ in } L^2((0, T); \mU). 
$$ 
Since 
$$
x_{T,\mathrm{opt}, k}(t) = e^{t A} \xi + J_k \int_0^t e^{(t-s)A} B u_{T,\mathrm{opt}, k}(s) \, ds,
$$
and
\begin{equation*}
\begin{split}
& (y_{T,\mathrm{opt}, k} - y_{T,\mathrm{opt}})'(t) = - A^* (y_{T,\mathrm{opt}, k}(t) - y_{T,\mathrm{opt}}(t))  - C^* C (x_{T,\mathrm{opt}, k}(t) - x_{T,\mathrm{opt}}(t)) \quad\mbox{ on } (0, T), \\
& (y_{T,\mathrm{opt}, k} - y_{T,\mathrm{opt}})(T) = 0, 
\end{split}
\end{equation*}
using the admissibility of $B$, we infer that 
$$
x_{T,\mathrm{opt}, k} \to x_{T,\mathrm{opt}} \ \mbox{ in } C([0, T]; \mH) \quad \mbox{ and } \quad y_{T,\mathrm{opt}, k} \to y \ \mbox{ in } C([0, T]; \mH).  
$$
This implies \eqref{lem-ap-cl2} and \eqref{lem-ap-cl3}. 
The proof is complete. 
\end{proof}

\subsection{Proof of \Cref{pro-technical}}
Given $\tau \in [0, T)$ and $\xi \in \mH$, we set
$$
\langle P_k (\tau) \xi, \xi \rangle  = \inf_{u \in L^2((\tau, T); \mU)} \int_0^T \left(|C x_k(t)|^2 + |u(t) |_{\mU}^2 \right)  dt + \langle P x_k(T), x_k(T) \rangle, 
$$
where $x_k(\cdot)$ is the unique solution of $x_k'(t) = A x_k(t) + B_k u(t)$ on $(\tau, T)$ such that $x_k(\tau) = \xi$ and $P$ is defined by \Cref{lem-P}. We set
$$
h_k(t) = y_{T,\mathrm{opt}, k}(t) - \bar y_k  -  P_{k}(t) (x_{T,\mathrm{opt}, k}(t) - \bar x_k) \quad \mbox{ and } \quad g_k(t)= h_k(T- t) \qquad\forall t\in [0, T].  
$$
Applying \Cref{lem1} and \Cref{lem-ap} with $z=0$ and $P_0= P$, we have 
$h_k \to h$ in $C([0, T]; \mH)$. 
This implies that
\be \label{pro-technical-pp-g}
g_k \to g \ \mbox{ in } C([0, T]; \mH). 
\ee
Given $\tau \in [0, T)$, we define $T_k(\tau)\in \cL(\mH)$ by $T_k(\tau) \xi = \varphi_k(\tau)$ where $\varphi_k'(t) = (-A^* + P_{k} B_k B_k^* ) \varphi_k(t)$ on $(\tau, T)$ with $\varphi_k(T) = \xi$ for every $\xi \in \mH$,
and we define $S_k(\tau) \in \cL(\mH)$ by $S_k(\tau) \eta = \psi_k(T)$ where $\psi_k'(t) = (A - B_k B_k^* P_{k}) \psi_k(t)$ on $(\tau, T)$ with $\psi_k(\tau) = \eta$ for every $\eta \in \mH$.

Applying \Cref{lem-transposition}, we have 
$$
\langle  \varphi_k(T), \psi_k(T) \rangle =  \langle \varphi(\tau), \psi(\tau) \rangle \qquad\forall \xi, \eta \in \mH, 
$$
hence
$$
\langle  T_k(\tau) \xi, \eta \rangle =  \langle \xi, S_k(\tau) \eta \rangle \qquad\forall \xi, \eta \in \mH.  
$$
Therefore,
\be \label{pro-technical-pp1}
T_k(\tau) = S_k(\tau)^*. 
\ee
Applying \Cref{lem-ap} in the time interval $[\tau, T]$ with $z=0$ and $P_0 = P$, we obtain
\be \label{pro-technical-pp2}
\lim_{k \to + \infty} S_k(\tau) \eta = S_{\infty,\mathrm{opt}}(T - \tau) \eta \ \mbox{ in $\mH$} \qquad\forall \eta \in \mH.  
\ee
Combining \eqref{pro-technical-pp1} and \eqref{pro-technical-pp2} yields that 
$T_k(\tau) \xi$ converges weakly to $S_{\infty,\mathrm{opt}}^*(T - \tau) \xi$ in $\mH$, for every $\xi \in \mH$.
Since, by \Cref{lem-technical2},
$h_k'(t)  = \left(-A^* + P_{k}B_k B_k^* \right) h_k(t)$ on $(0, T)$,
we infer from \eqref{pro-technical-pp-g} that 
$g(t) = S_{\infty,\mathrm{opt}}^*(t) g(0)$ for every $t\in[0, T]$.
The proof is complete.

\section{Proof of \Cref{thm-main}} \label{sect-main}
In what follows in this proof, $c$ is a generic positive constant depending only on $(A, B, C)$ (and on $T_0$, which depends on $(A, B, C)$ as well) and which can change from one place to another. 

First, recalling that $(x_{T,\mathrm{opt}}-\bar x,y_{T,\mathrm{opt}}-\bar y)$ satisfies the system \eqref{lem-technical-Sys-mod} and that $x_{T,\mathrm{opt}}(0)=x_0$ and $y_{T,\mathrm{opt}}(T)=0$, applying \Cref{lem-transposition}, we have 
\begin{equation*}
\begin{split}
&\ \langle x_0 - \bar x, y_{T,\mathrm{opt}}(0) - \bar y \rangle + \langle x_{T,\mathrm{opt}}(T) - \bar x, \bar y \rangle \\
=&\ \langle x_{T,\mathrm{opt}}(0) - \bar x, y_{T,\mathrm{opt}}(0) - \bar y \rangle - \langle x_{T,\mathrm{opt}}(T) - \bar x, y_{T,\mathrm{opt}}(T) - \bar y \rangle \\
=&\ \int_0^T \left( |B^* (y_{T,\mathrm{opt}}(t) - \bar y)|^2 +  |C (x_{T,\mathrm{opt}}(t) -  \bar x)|^2 \right) dt, \\  
=&\ \int_0^T \left( | u_{T,\mathrm{opt}}(t) - \bar u|^2 +  |C (x_{T,\mathrm{opt}}(t) -  \bar x)|^2 \right) dt ,
\end{split}
\end{equation*}
from which it follows, by the Cauchy-Schwarz inequality, that
\begin{equation}\label{thm-main-p1}
\int_0^T \left( | u_{T,\mathrm{opt}}(t) - \bar u|^2 +  |C (x_{T,\mathrm{opt}}(t) -  \bar x)|^2 \right) dt  \leq \vert x_0-\bar x\vert \vert y_{T,\mathrm{opt}}(0) - \bar y\vert + \vert x_{T,\mathrm{opt}}(T) - \bar x\vert \vert \bar y\vert .
\end{equation}

Second, since $(A^*, B^*)$ is finite-time observable in time $T_0$, using the arguments in \cite[Remark 2.1]{PZ13}, we infer from the equation of $y_{T,\mathrm{opt}} - \bar y$ in \eqref{lem-technical-Sys-mod} that, for any $T \ge T_0$, 
\begin{equation} \label{thm-main-p2}
\begin{split}
|y_{T,\mathrm{opt}}(0) - \bar y |^2 
&\le c \int_0^{T_0} \left( | B^*(y_{T,\mathrm{opt}}(t) - \bar y)|^2 + |C (x_{T,\mathrm{opt}}(t) -  \bar x)|^2 \right) dt \\
&= c \int_0^{T_0} \left( | u_{T,\mathrm{opt}}(t) - \bar u |^2 +  |C (x_{T,\mathrm{opt}}(t) -  \bar x)|^2 \right) dt \\
&\le  c \int_0^T \left( | u_{T,\mathrm{opt}}(t) - \bar u |^2 +  |C (x_{T,\mathrm{opt}}(t) -  \bar x)|^2 \right) dt .
\end{split}
\end{equation}

Third, similarly, since $(A, C)$ is finite-time observable in time $T_0$, using the arguments in \cite[Remark 2.1]{PZ13}, we infer from the equation of $x_{T,\mathrm{opt}} - \bar x$ in \eqref{lem-technical-Sys-mod} that, for any $T \ge T_0$, 
\begin{equation} \label{thm-main-p3}
\begin{split}
|x_{T,\mathrm{opt}}(T) - \bar x |^2 
&\le c \int_{T-T_0}^T \left( | C(x_{T,\mathrm{opt}}(t) - \bar x)|^2 +  |u_{T,\mathrm{opt}}(t) -  \bar u|^2  \right) dt \\
&\le c \int_0^T \left( | C(x_{T,\mathrm{opt}}(t) - \bar x)|^2 +  |u_{T,\mathrm{opt}}(t) -  \bar u|^2  \right) dt .
\end{split}
\end{equation}
Combining \eqref{thm-main-p1}, \eqref{thm-main-p2} and \eqref{thm-main-p3} yields 
$$
\int_0^T \left( | u_{T,\mathrm{opt}}(t) - \bar u|^2 +  |C (x_{T,\mathrm{opt}}(t) -  \bar x)|^2 \right) dt 
\le c \left (|x_0-\bar x|^2 + |\bar y|^2 \right)  .    
$$
This in turn implies, by \eqref{thm-main-p2} and \eqref{thm-main-p3}, that
\be \label{thm-main-p3-***}
| x_{T,\mathrm{opt}} (T) -  \bar x | \le c \left (|x_0-\bar x| + |\bar y| \right) \quad \mbox{ and } \quad |y_{T,\mathrm{opt}}(0) -  \bar y| \le c \left (|x_0-\bar x| + |\bar y| \right) . 
\ee
Recalling that $h$ and $g$ are defined by \eqref{pro-technical-def-hg}, by \Cref{pro-technical}, we have  
$g(t) = S_{\infty,\mathrm{opt}}^* (t) g(0)$ for every $t\ge 0$, where $S_{\infty,\mathrm{opt}}$ is the semigroup given in \Cref{pro-def-Sopt}. 
As a consequence, we infer from \Cref{pro-def-Sopt} that  
\begin{equation} \label{thm-main-decay1}
\begin{split}
|y_{T,\mathrm{opt}}(t) - \bar y  -  P (x_{T,\mathrm{opt}}(t) - \bar x)|
&\le c e^{-\lambda (T - t)} |y_{T,\mathrm{opt}}(T) - \bar y  -  P (x_{T,\mathrm{opt}}(T) - \bar x)| \\
&\mathop{\le}^{\eqref{thm-main-p3-***}} c e^{-\lambda (T - t)} \left (|x_0-\bar x| + |\bar y| \right)
\end{split}
\end{equation}
for every $t\in [0, T]$.

Setting
\be \label{thm-main-decay-def-xieta}
\xi(t) = x_{T,\mathrm{opt}}(t) - \bar x \quad \mbox{ and } \quad \eta(t) = y_{T,\mathrm{opt}}(t) - \bar y \qquad\forall t\in[0, T],
\ee
we have, on $(0,T)$, 
\begin{equation*} 
\begin{split}
\xi'(t) &= A \xi(t) - B B^* \eta(t), \\
\eta'(t) &= - A^* \eta(t) - C^* C \xi(t).
\end{split}
\end{equation*}
Applying \Cref{lem-transposition}, we obtain, for any $\tau \in [0, T]$,  
$$
\langle \xi(\tau), \eta (\tau) \rangle - \langle \xi(0), \eta (0) \rangle
= - \int_0^{\tau} \left(|B^* \eta(t)|^2 + |C \xi(t)|^2 \right) dt. 
$$
Using that, by \eqref{thm-main-decay1}, $|\eta(t) - P \xi(t) | \mathop{\le}^{} c e^{- \lambda (T - t)} \left (|x_0-\bar x| + |\bar y| \right)$ for every $t \in [0, T]$,
we get 
\begin{multline}\label{thm-main-decay2-1}
\langle \xi(\tau), P \xi (\tau) \rangle +  \int_0^{\tau} \left(|B^* \eta(t)|^2 + |C \xi(t)|^2 \right) dt \\
\le \langle \xi(0), P \xi (0) \rangle  + c e^{- \lambda (T - \tau)} \left (|x_0-\bar x| + |\bar y| \right) |\xi (\tau)| \qquad\forall \tau \in [0, T]. 
\end{multline}

Let $\hu_{\infty,\mathrm{opt}, \xi(0)}$ (resp., $\hu_{\infty,\mathrm{opt}, \xi(\tau)}$) be the optimal control solution of 
$$
\inf_{u \in L^2((0, + \infty); \mU)}  \int_0^\infty \left(|Cx(t)|^2 + |u(t) |^2 \right) dt
$$
where $x(\cdot)$ is the unique solution of $x'(t) = Ax(t) + Bu(t)$ on $(0,+\infty)$ such that $x(0) = \xi(0)$ (resp., such that $x(0) = \xi(\tau)$), and let $\hx_{\infty,\mathrm{opt}, \xi(0)}$ (resp., $\hx_{\infty,\mathrm{opt}, \xi(\tau)}$) be the corresponding solution.  

On the one part, by definition of $P$, we have
\begin{multline} \label{thm-main-decay2-4}
\langle \xi(0), P \xi(0) \rangle  =  \int_0^{\infty} \left(|\hu_{\infty,\mathrm{opt}, \xi(0)}(t)|^2 + |C \hx_{\infty,\mathrm{opt}, \xi(0)}(t)|^2 \right) dt \\
= \inf_{u \in L^2((0, + \infty); \mU)}  \int_0^\infty \left(|Cx(t) |^2 + |u(t) |^2 \right) dt .
\end{multline}

On the other part, we define a new trajectory, by concatenation: 
\be
\xi_{e, \tau} (t) =  \left\{\begin{array}{ll} 
\xi (t) & \mbox{ in } (0, \tau), \\
\hx_{\infty,\mathrm{opt}, \xi(\tau)} (t - \tau) & \mbox{ in } (\tau, + \infty), 
\end{array} \right. \quad \mbox{ and } \quad u_{e, \tau}(t) =  \left\{\begin{array}{ll} 
- B^*\eta (t) & \mbox{ in } (0, \tau), \\
\hu_{\infty,\mathrm{opt}, \xi(\tau)} (t - \tau) & \mbox{ in } (\tau, + \infty),
\end{array} \right. 
\ee
so that
\be \label{thm-main-decay2-2}
\xi_{e, \tau}'(t) = A \xi_{e, \tau}(t) + B u_{e, \tau}(t)\quad\mbox{on } (0, + \infty), \qquad
\xi_{e, \tau} (0) = \xi(0),  
\ee
and, by \eqref{thm-main-decay2-1}, we have
\be \label{thm-main-decay2-3}
\int_0^{\infty} \left(|u_{e, \tau}(t)|^2 + |C \xi_{e, \tau}(t)|^2 \right) dt  \le \langle \xi(0), P \xi (0) \rangle + c e^{- \lambda (T - \tau)} \left (|x_0-\bar x| + |\bar y| \right) |\xi (\tau)|. 
\ee

We now use again the parallelogram identity \eqref{thm-main-Parallelogram} (as in the proof of Lemma \ref{lem-ap}).
We infer from \eqref{thm-main-decay2-4}, \eqref{thm-main-decay2-2}, and \eqref{thm-main-decay2-3} that 
$$
\int_0^{\infty} \left(\left|u_{e, \tau}(t) - u_{\infty,\mathrm{opt}, \xi(0)}(t)\right|^2 + \left|C ( \xi_{e, \tau}(t) - \hx_{\infty,\mathrm{opt}, \xi(0)}(t)) \right|^2\right) dt \le c e^{- \lambda (T - \tau)} \left (|x_0-\bar x| + |\bar y| \right) |\xi (\tau)|. 
$$
This implies, in particular, that
\begin{multline} \label{thm-main-decay2-5}
\int_0^{\tau} \left(\left|(u_{T,\mathrm{opt}}(t) - \bar u) - u_{\infty,\mathrm{opt}, \xi(0)}(t)\right|^2 + \left|C \left( (x_{T,\mathrm{opt}}(t) - \bar x) - \hx_{\infty,\mathrm{opt}, \xi(0)}(t) \right) \right|^2 \right) dt \\ 
\le c e^{- \lambda (T - \tau)}\left (|x_0-\bar x| + |\bar y| \right) |\xi (\tau)|. 
\end{multline}
Since, for $\tau \ge T_0$,  
$$
\int_{\tau - T_0}^{\tau} \left(\left|u_{\infty,\mathrm{opt}, \xi(0)}(t)\right|^2 + \left|C  \hx_{\infty,\mathrm{opt}, \xi(0)}(t)\right|^2 \right) dt  \le c e^{- 2 \lambda \tau} |\xi(0)|^2,  
$$
we infer from \eqref{thm-main-decay2-5} that 
\begin{multline}\label{thm-main-decay2-7}
\int_{\tau - T_0}^{\tau} \left(\left|u_{T,\mathrm{opt}}(t) - \bar u\right|^2 + \left|C (x_{T,\mathrm{opt}}(t) - 
\bar x)\right|^2 \right) dt  \\
\le c \left(e^{- \lambda (T - \tau)} |\xi(\tau)| + e^{-2\lambda \tau} \left (|x_0-\bar x| + |\bar y| \right) \right) \left (|x_0-\bar x| + |\bar y| \right). 
\end{multline}
Since  $(A, C)$ is finite-time observable in time $T_0$, we infer from the equation of $x_{T,\mathrm{opt}} - \bar x$ in \eqref{lem-technical-Sys-mod} that, for  $T_0 < \tau  \le T$, 
\be \label{thm-main-decay2-8}
|x_{T,\mathrm{opt}}(\tau) - \bar x |^2 \le c \int_{\tau-T_0}^\tau \left( | C(x_{T,\mathrm{opt}}(t) - \bar x)|^2 +  |u_{T,\mathrm{opt}}(t) -  \bar u|^2 \right) dt.  
\ee
Combining \eqref{thm-main-decay2-7} and \eqref{thm-main-decay2-8} gives
\be \label{thm-main-decay2}
|x_{T,\mathrm{opt}}(\tau) - \bar x | \le c \big(e^{- \lambda (T - \tau)} + e^{-\lambda \tau} \big) \left (|x_0-\bar x| + |\bar y| \right) \qquad\forall\tau\in[T_0,T].
\ee
It is clear that \eqref{thm-main-decay2} holds for $0 \le \tau \le T_0$. 
Combining \eqref{thm-main-decay1} and \eqref{thm-main-decay2} yields 
\be \label{thm-main-decay3}
|x_{T,\mathrm{opt}}(\tau) - \bar x | + |y_{T,\mathrm{opt}}(\tau) - \bar y | \le c \big(e^{- \lambda (T - \tau)} + e^{-\lambda \tau} \big) \left (|x_0-\bar x| + |\bar y| \right)  \qquad\forall\tau\in[0,T].
\ee

We finally estimate $u_{T,\mathrm{opt}} - \bar u$. Let $t \in [0, T/2]$. 
Similarly to what has been done at the beginning of this section, applying \Cref{lem-transposition}, we have 
\begin{multline*}
\langle x_{T,\mathrm{opt}}(t) - \bar x, y_{T,\mathrm{opt}}(t) - \bar y \rangle - \langle x_{T,\mathrm{opt}}(T-t) - \bar x, y_{T,\mathrm{opt}}(T-t) - \bar y \rangle \\
= \int_t^{T-t} \left( | u_{T,\mathrm{opt}}(t) - \bar u|^2 +  |C (x_{T,\mathrm{opt}}(t) -  \bar x)|^2 \right) dt ,
\end{multline*}
from which it follows, by the Cauchy-Schwarz inequality, that
$$
\int_t^{T-t} | u_{T,\mathrm{opt}}(t) - \bar u|^2 \, dt  \leq 
\vert x_{T,\mathrm{opt}}(t) - \bar x\vert \vert y_{T,\mathrm{opt}}(t) - \bar y \vert + \vert x_{T,\mathrm{opt}}(T-t) - \bar x \vert \vert y_{T,\mathrm{opt}}(T-t) - \bar y \vert .
$$
Using \eqref{thm-main-decay3}, the desired estimate follows.
The proof is complete. \qed

\appendix

\section{Proof of \Cref{pro-Opt-T}} \label{sect-pro-Opt-T}

  \setcounter{equation}{0}  
  \setcounter{lemma}{0}  

Before proving \Cref{pro-Opt-T}, we establish a preliminary result. 

\begin{lemma}\label{lem-Opt-T} 
Let $T>0$. Given $x_0, z  \in \mH$, consider the pair $(\tx, \tu)$ where $\tu \in L^2((0, T); \mU)$ and $\tx$ is the unique solution of $\tx'(t) = A \tx(t) + B \tu(t)$ on $(0, T)$ such that $\tx(0) = x_0$. 
Then $(\tx, \tu)$ is the optimal solution of \eqref{optimal-control-T} if and only if 
\be \label{lem-Opt-T-tu}
\tu= - B^* \ty\ \mbox{ in } L^2((0, T); \mU),  
\ee
where $\ty$ is the unique solution of $\ty'(t)=  - A^* \ty(t) - C^* (C \tx(t) - z)$ on $(0, T)$ such that $\ty (T) = P_0 \tx(T)$. 
\end{lemma}

\begin{remark} \rm The equality \eqref{lem-Opt-T-tu} can written as
\be \label{rem-lem-Opt-T-u}
\tu (t)= - B^* \left(e^{(T-t)A^*} P_0 \tx(T) + \int_t^T e^{(s-t) A^*} C^* (C \tx(s) - z) \, ds\right) \quad\mbox{ on } (0, T).
\ee
It is a consequence of \Cref{lem-transposition} that the left-hand side of \eqref{rem-lem-Opt-T-u} belongs to $L^2((0, T); \mU)$ as mentioned in \Cref{rem-B*-L2}.  
\end{remark}

\begin{proof}[Proof of \Cref{lem-Opt-T}] Let $u \in L^2((0, T); \mU)$ and let $x(\cdot)$ be the unique solution of 
$x'(t) = A x(t) + B u(t)$ on $(0, T)$ such that $x(0) = x_0$.
We have 
\begin{multline*} 
\int_0^T \left( |Cx(t) - z|^2 + |u(t)|^2  \right) dt - \int_0^T \left( |C \tx(t) - z|^2 + |\tu(t)|^2  \right) dt \\
= \int_0^T \left( |C (x(t) - \tx(t))|^2 + |u(t) - \tu|^2  \right) dt + \langle P_0 (x(T) - \tx(T)), x(T) - \tx(T) \rangle \\
+ 2 \int_0^T \left( \langle C^* (C \tx(t) - z),  (x(t) - \tx(t)) \rangle  +\langle  \tu(t), u(t) - \tu(t) \rangle  \right) \, dt  + 2 \langle P_0 \tx(T), x (T) - \tx(T) \rangle. 
\end{multline*}
Hence, $(\tx, \tu)$ is the optimal solution if and only if 
\be \label{thm-Opt-T-p1}
\int_0^T \left( \langle C^* (C \tx(t) - z),  Lv  \rangle  +\langle  \tu(t), v \rangle \right) dt  +  \langle P_0 \tx(T), L v(T) \rangle = 0 
\qquad\forall v \in L^2((0, T); \mU), 
\ee
where, for every $v \in L^2((0, T); \mU)$, we denote by $Lv$ the solution of 
$(Lv)'(t) = A (L v)(t) + B v(t)$ on $(0, T)$ such that $(Lv)(0) = 0$.

Applying \Cref{lem-transposition} to $Lv$ and $\ty$, we have 
$$
\langle Lv(T), \ty (T) \rangle - \langle Lv(0), \ty (0) \rangle = \int_0^T \langle v(t), B^* \ty(t) \rangle \, dt - \int_0^T \langle Lv(t), C^* (C \tx(t) - z) \rangle \, dt.  
$$
Since $Lv(0) = 0$ and $\ty(T) = P_0 \tx(T)$, it follows that  
\be \label{thm-Opt-T-p2}
\langle Lv(T), P_0 \tx(T) \rangle   +  \int_0^T \langle Lv(t), C^* (C \tx(t) - z) \rangle \, dt = \int_0^T \langle v(t), B^* \ty(t) \rangle \, dt.   
\ee
Combining \eqref{thm-Opt-T-p1} and \eqref{thm-Opt-T-p2} yields 
$$
\int_0^T \langle \tu(t) + B^* \ty(t), v(t) \rangle \, dt = 0 \qquad\forall v \in L^2((0, T); \mU). 
$$
The conclusion follows. 
\end{proof}

We are now in a position to prove \Cref{pro-Opt-T}. 

\begin{proof}[Proof of \Cref{pro-Opt-T}]
By \Cref{lem-Opt-T}, it remains to prove that when $z =0$, we have $\ty(t) = P_T(t) \tx (t)$ for every $t\in [0, T]$, where $P_T(t)$ is defined by \eqref{def-Pt}. 

Applying \Cref{lem-transposition} to $\tx$ and $\ty$ on the time interval $(\tau, T)$, we have 
\be \label{pro-Opt-T-1-1}
\langle \tx(T), \ty (T) \rangle  - \langle \tx(\tau), \ty (\tau) \rangle = \int_{\tau}^{T} \left( \langle \tu(t), B^*\ty(t) \rangle - \langle \tx(t), C^*C \tx(t) \rangle \right) dt. 
\ee
Since $\tu = - B^* \ty$ on $(0, T)$ by \eqref{pro-Opt-T-sys-opt2} and $\ty(T) = P_0 \tx(T)$, we infer from \eqref{pro-Opt-T-1-1} that 
$$
\langle \tx(\tau), \ty (\tau) \rangle =   \int_{\tau}^{T } \left(|\tu(t)|^2 + |C \tx(t)|^2 \right) dt + \langle P_0 \tx(T), \tx(T) \rangle \qquad\forall \tau \in [0, T). 
$$
By definition of $P_T(\tau)$, we thus have 
\be \label{pro-Opt-T-tx-ty1}
\langle \tx(\tau), \ty (\tau) \rangle =   \langle P_T(\tau) \tx(\tau), \tx(\tau) \rangle \qquad\forall\tau\in [0, T). 
\ee
In particular, $\langle \tx(0), \ty (0) \rangle =   \langle P_T(0) \tx(0), \tx(0) \rangle$.

One can check, by uniqueness of the optimal control, that $\ty(0)$ is a linear function of $\tx(0)$, and moreover this linear function is continuous.  It follows that there exists $M_0 \in \cL(\mH)$ such that $\ty(0) = M_0 \tx(0)$. 

We claim that $M_0$ is symmetric and thus obtain that
$\ty(0) = M_0 \tx(0)$.
Indeed, let $\tu_1$ and $\tu_2$ be two optimal controls corresponding to the initial data $\xi_1$ and $\xi_2$ at  time $t = 0$, respectively, and let $\tx_1$ and $\tx_2$ be the corresponding solutions and $\ty_1$ and $\ty_2$ be the corresponding adjoint states. 
Applying \Cref{lem-transposition}, we have
\be \label{pro-Opt-T-m1}
\langle \tx_1(T), \ty_2 (T) \rangle - \langle \tx_1(0), \ty_2 (0) \rangle = - \int_0^T \left( \langle B^* \ty_1(t), B^* \ty_2(t) \rangle + \langle C \tx_1(t),  C \tx_2(t) \rangle \right) dt
\ee
and 
\be \label{pro-Opt-T-m2}
\langle \tx_2(T), \ty_1 (T) \rangle - \langle \tx_2(0), \ty_1 (0) \rangle = - \int_0^T \left( \langle B^*\ty_2(t), B^*\ty_1(t) \rangle + \langle C \tx_2(t),  C \tx_1(t) \rangle \right) dt.  
\ee
Since, by \eqref{pro-Opt-T-sys-opt1} 
$$
\langle \tx_1(T), \ty_2 (T) \rangle = \langle \tx_1(T), P_0 \tx_2 (T) \rangle, \quad  \langle \tx_2(T), \ty_1 (T) \rangle = \langle \tx_2(T), P_0 \tx_1 (T) \rangle, 
$$ 
and $P_0$ is symmetric, 
we infer from \eqref{pro-Opt-T-m1} and \eqref{pro-Opt-T-m2} that 
$\langle \tx_1(0), \ty_2 (0) \rangle = \langle \tx_2(0), \ty_1 (0) \rangle$.
Hence $\langle \xi_1, M_0 \xi_2 \rangle = \langle \xi_2, M_0 \xi_1 \rangle$.
Since $\xi_1, \xi_2\in\mH$ are arbitrary, we conclude that $M_0 \in \cL(\mH)$ is symmetric, and the claim is proved. 

By changing the starting time, we have thus proved that,  for every $\tau \in [0, T)$, 
\be \label{pro-Opt-T-tx-ty2}
\ty(\tau) = M(\tau) \tx(\tau)\quad \mbox{ for some symmetric } M(\tau) \in \cL(\mH).  
\ee
As a consequence of \eqref{pro-Opt-T-tx-ty1} and \eqref{pro-Opt-T-tx-ty2}, we infer that $M(\tau) = P_T(\tau)$ in $[0, T)$, which yields $\ty(\tau) = P_T(\tau) \tx (\tau)$ for every $\tau\in[0, T]$.
The proof is complete. 
\end{proof}

\providecommand{\bysame}{\leavevmode\hbox to3em{\hrulefill}\thinspace}
\providecommand{\MR}{\relax\ifhmode\unskip\space\fi MR }
\providecommand{\MRhref}[2]{%
  \href{http://www.ams.org/mathscinet-getitem?mr=#1}{#2}
}
\providecommand{\href}[2]{#2}

\end{document}